\documentclass[journal,twoside,web]{ieeecolor}
\usepackage{generic}
\usepackage{cite}
\usepackage{amsmath,amssymb,amsfonts}
\usepackage{algorithmic}
\usepackage{graphicx}
\usepackage{algorithm,algorithmic}
\usepackage{hyperref}
\hypersetup{hidelinks=true}
\usepackage{textcomp}
\usepackage{tikz}
\usepackage{epstopdf}
\usepackage{stfloats}

\newtheorem{theorem}{Theorem}[section]

\newtheorem{lemma}[theorem]{Lemma}
\newtheorem{proposition}[theorem]{Proposition}

\newtheorem{remark}[theorem]{Remark}
\newtheorem{assumption}[theorem]{Assumption}

\DeclareMathOperator*{\esssup}{ess\,sup}

\def\BibTeX{{\rm B\kern-.05em{\sc i\kern-.025em b}\kern-.08em
    T\kern-.1667em\lower.7ex\hbox{E}\kern-.125emX}}
\begin{document}
\title{Micro-Macro Backstepping Control of Large-Scale Hyperbolic Systems$^*$}
\author{Jukka-Pekka Humaloja, \IEEEmembership{Senior Member, IEEE}, and Nikolaos 
Bekiaris-Liberis, 
\IEEEmembership{Senior Member, IEEE}
\thanks{$^*$Funded by the European Union (ERC, C-NORA, 101088147). Views and opinions 
	expressed are however those of the authors only and do not necessarily reflect those of the 
	European Union or the European Research Council Executive Agency. Neither the European 
	Union nor the granting authority can be held responsible for them.}
\thanks{The authors are with the Department of Electrical and Computer Engineering, 
	Technical University of Crete, Chania, Greece. Emails: jhumaloja@tuc.gr and nlimperis@tuc.gr.}
}

\maketitle

\begin{abstract}
We introduce a control design and analysis framework for micro-macro, boundary control of 
large-scale, $n+m$ hyperbolic PDE systems. Specifically, we develop feedback laws for 
stabilization of hyperbolic systems at the micro level (i.e., of the large-scale system) that employ 
a) measurements obtained from the $n+m$ system (i.e., at micro level) and kernels constructed 
based on an $\infty+\infty$ continuum system counterpart (i.e., at macro level), or b) kernels and 
measurements both stemming from a continuum counterpart, or c) averaged-continuum 
kernels/measurements. We also address (d)) stabilization of the continuum (macro) system, 
employing continuum kernels and measurements. The significance of addressing a)--d) lies in the 
facts that for large-scale hyperbolic systems computation of stabilizing control kernels 
(constructed for the n+m system) may become intractable and in different applications only 
average (macro) measurements may be available. The main design and analysis steps involved in 
a)--d) are the following. Towards addressing d) we derive in a constructive manner an 
$\infty+\infty$ continuum approximation of $n+m$ hyperbolic systems and establish that its 
solutions approximate, for large $n$ and $m$, the solutions of the $n+m$ system. We then 
construct a feedback law for stabilization of the $\infty+\infty$ system via introduction of a 
continuum-PDE backstepping transformation. We establish well-posedness of the resulting 4-D 
kernel equations and prove closed-loop stability via construction of a novel Lyapunov functional. 
Furthermore, under control configuration a) we establish that the closed-loop system is 
exponentially stable provided that $n$ and $m$ are large, by proving that the exact, stabilizing 
$n+m$ control kernels can be accurately approximated by the continuum kernels. While under 
control configurations b) and c), we establish closed-loop stability capitalizing on the established 
solutions' and kernels' approximation properties via employment of infinite-dimensional ISS 
arguments. We provide two numerical simulation examples to illustrate the effectiveness and 
potential limitations of our design approach.
\end{abstract}


\section{Introduction} \label{sec:intro}

\subsection{Motivation}

\begin{table*}[tp] 
\begin{center}
\begin{tabular}{|c|c|c|c|c|} \hline & & & &  \\ [-8pt]
 & Objective in & Level of control & Kernels & Measurements	\\ 
& present paper & implementation & construction & available  \\  \hline & & & &  \\ [-8pt]
New result & Control of continuum & Macro & Macro & Macro \\ 
& hyperbolic PDEs ($n=m = \infty$) & & & \\ \hline & & & &  \\ [-8pt]
New result & Control of large-scale & Micro & Macro & Micro \\ 
& hyperbolic PDEs (finite/large $n,m$) & & & \\ \hline & & & &  \\ [-8pt]
New result & Control of large-scale & Micro & Micro & Macro \\ 
& hyperbolic PDEs (finite/large $n,m$) & & & \\ \hline & & & &  \\ [-8pt]
New result & Control of large-scale & Micro & Averaged macro & Averaged macro \\ 
& hyperbolic PDEs (finite/large $n,m$) & & &  \\ \hline & & & &  \\ [-8pt]
Existing literature & Control of hyperbolic & Micro & Micro & Micro \\ 
\cite{AnfAamBook, AurBre22, HuLDiM16, AurDiM16, CorHuL17, DiMArg18, HuLVaz19, CorHuL21, 
GehDeu25, RamZwaCDC17} 
&  PDEs 
(finite
$n,m$) 
& & & \\ \hline
\end{tabular}
\end{center}
\caption{Overview of Problems Addressed and of Contributions}
\label{tab:ccc}
\end{table*}

\IEEEPARstart{M}{icro-macro} control, i.e., the approach in which control is implemented at or 
designed for different system levels, or employs measurements stemming from different system 
levels, has been, heretofore, introduced only for specific applications, in particular, traffic flow 
control and only for the cases where the underlying models considered consist of large-scale ODE 
systems; see, for example, \cite{GoaRos17, KarThe22, SpiMan18, BekDel21, MirIov23}. Taking a 
significant step forward, in this paper, we introduce a new and systematic framework for design 
and analysis of micro-macro controllers for large-scale hyperbolic PDE systems. 

The reasons that 
such a general (i.e., not developed only for a specific engineering application) approach for 
micro-macro control of large-scale PDE systems is significant stem from the facts that such a 
setup appears in different applications and it enables introduction of new control design and 
analysis ideas/tools. In particular, such an approach may be essential when dealing with 
large-scale hyperbolic PDE systems in order to construct feedback laws that are computationally 
tractable and that rely on availability of only some average (macro) measurements. Among other 
applications, such an approach may be suitable for lane-free traffic \cite{MalPap21} or 
large-scale traffic networks \cite{HerKla03, YuHKrs21, BurYuH21}, large-scale blood 
flow networks \cite{BikPhd, ReyMer09}, and large-scale epidemics spreading networks 
\cite{GuaPri20, KitBes22}; all of which can be described by systems consisting of large-scale (or 
continua-of) hyperbolic PDEs and for which control design/implementation or measurements may 
be available at different system levels (micro or macro). 

 To help the reader better understand the setup of each problem we address, as shown in 
 Table~\ref{tab:ccc}, we explain 
 how these may emerge in traffic flow control-related applications. The 
 setup in the first row of Table~\ref{tab:ccc} may emerge, for example, in the case of lane-free 
 traffic flow, 
 where traffic is viewed/modeled as a 2-D continuum/fluid \cite{MalPap21}. The setup 
 in which 
 for a given  large-scale system only macro measurements are available (corresponding to the 
 problem in the third  row of Table~\ref{tab:ccc}), may appear, for example, in cases when control 
 is performed via  manipulation of
individual, automated/connected vehicles’ trajectories, based 
on density (or speed) 
measurements or estimates (in a given road segment) that correspond to some average
spacing between individual vehicles (or to some average speed of vehicles). Moreover, the setup 
where only averaged macro measurements (fourth row of Table~\ref{tab:ccc}) are available may 
be motivated by, for example, large-scale transportation networks where only averaged (over a 
given network segment) measurements of density/speed may be available, whereas a
traffic controller is implemented locally, at each traffic system component, see, e.g., 
\cite{TumCan22, SirGer18}. 

\subsection{Literature}

Even though there is no approach specifically addressing the problem of micro-macro control of 
large-scale and continua-of hyperbolic PDEs, the most closely related literature includes the 
results on control of specific classes of large-scale hyperbolic systems via a continuum PDE 
approach utilizing backstepping \cite{HumBek25b, HumBek26} and the results in \cite{AllKrs25} 
dealing with backstepping control of a specific continuum of hyperbolic PDEs. In particular, in the 
former results, control of $n+1$ and $n+m$ (for large $n$ and $m$) hyperbolic systems is 
considered, together with control of their continuum $\infty+1$ and $\infty+m$ counterparts, 
respectively; whereas in the latter results, the case of $\infty+1$ continua is addressed. As we use 
the backstepping design concept, the results in \cite{AnfAamBook, AurBre22, HuLDiM16, 
AurDiM16, CorHuL17, DiMArg18, HuLVaz19, CorHuL21,GehDeu25, RamZwaCDC17} concerning 
backstepping-based control 
of $n+m$ hyperbolic PDEs are also relevant, even though they do not specifically address 
large-scale or continua-of hyperbolic PDE systems (or the exact interplay between them). 

Furthermore, since a main motivation for micro-macro, PDE backstepping control design is that it 
enables construction/computation of stabilizing, backstepping kernels for large-scale PDE 
systems, with computational complexity that does not grow with the number of PDE states 
components, the results in \cite{BhaShi24, QiJZha24, WanDia25} concerning computation of 
backstepping kernels via neural operators for single/two-component PDE systems; the results in 
\cite{AurMor19} that present a late-lumping-based approach; and the results in 
\cite{VazCheCDC23} that rely on power series representations for computation of the kernels, are 
also relevant. We note that these results do not address micro-macro control and do not aim at 
addressing the growing computational complexity of backstepping kernels as the number of PDE 
states components becomes large (and thus, computational complexity in these approaches may 
still grow with the number of system components). Finally, although the technical tools we 
develop and utilize here are different, we also consider as relevant results dealing with control of 
large-scale ODE systems via a continuum approach, such as, for example, \cite{FerBuf06, 
NikCan22, MeuKrs11, QiJVaz15, FriKrs11, ZhaVaz24}, as we borrow the idea of constructing PDE 
continua for design of controllers for the original, large-scale PDE systems considered.

\subsection{Contributions}

\paragraph{Conceptual contributions} In this paper we develop a new control design approach 
for micro-macro control of large-scale hyperbolic systems. In the framework we introduce there 
are different configurations, which we specify, for micro-macro control depending on which level 
(micro or macro) control is applied, essentially corresponding to whether the objective is control 
of the macro (continuum) or micro (large-scale) system; on which level measurements are 
obtained, i.e., on whether measurements are available directly from the large-scale (micro) 
system or they are available only on average (based on a macro system counterpart); and on 
which level control gains are constructed, i.e., whether control gains are constructed based on the 
continuum (macro) system or based on the large-scale (micro) system. In the present paper we 
address the control design and analysis problems shown in Table~\ref{tab:ccc}, namely, we 
consider the case where the objective is stabilization of (micro) large-scale, $n+m$, linear 
hyperbolic systems utilizing 
control kernels that are constructed on a macro level (i.e., based on a continuum, $\infty+\infty$
system) and/or measurements that are obtained on a macro level (i.e., from a continuum 
$\infty+\infty$ system counterpart) or even as averaged macro measurements. We also address 
the problem of stabilization of the continuum (macro) system itself, using continuum (macro) 
kernels and measurements. 

\paragraph{Technical contributions} To execute the above conceptual ideas we have to 
introduce a new control design and analysis approach, whose main ingredients are the following. 
The first step for computing stabilizing continuum kernels in a computational tractable manner, 
i.e., based on a continuum system, is to actually construct/derive a proper continuum system that 
approximates (in certain sense) the original large-scale system. We resolve this problem by 
introducing a constructive approach for constructing continuum systems based on a given 
large-scale, $n+m$ hyperbolic system, as $n$ and $m$ tend to infinity. We then establish that 
the solutions of 
the continuum $\infty+\infty$ system approximate (in a certain sense) the solutions of the 
large-scale $n+m$ system 
for large $n$ and $m$. The second step/contribution is to develop a backstepping state-feedback 
control law to exponentially stabilize the class of $\infty+\infty$ continua of hyperbolic systems. 
We achieve this via introduction of a continuum-PDE infinite-dimensional backstepping 
transformation. The key technical challenges that we resolve in our approach are the study of 
well-posedness of the resulting kernel equations and the construction of a novel Lyapunov 
functional, as neither of them follows from the existing results (although we rely on specific 
existing results, in particular, on \cite{HumBek25b, HumBek26}, as reasonably expected). In 
particular, well-posedness of 
the kernel equations does not follow in a straightforward manner from existing results as the 
backstepping procedure we consider gives rise to two continuum kernel equations evolving on a 
4-D domain, which is obtained by a continuation of the prismatic 2-D domain of $n+m$ kernel 
equations over a 2-D function space in $L^2([0; 1]^2;\mathbb{R})$. We note that both of the 
above results 
constitute a significant step forward as compared with the case of the respective construction for 
only large $n$ from \cite{HumBek26}, as the case where $m\to\infty$ imposes unique technical 
challenges, due 
to the input space becoming infinite-dimensional, which mainly arise because the pointwise 
arguments employed in $\mathbb{R}^m$ are not viable for $L^2$ functions and the kernel 
equations evolve on a 4-D domain.

We next provide two more key results. In the first, we develop a micro-macro control design 
methodology for 
stabilization of the large-scale, $n+m$ hyperbolic system utilizing micro measurements and 
macro kernels, i.e., kernels constructed based on the continuum $\infty+\infty$ PDE system. We 
establish that the closed-loop system is exponentially stable provided that $n$ and $m$ are 
sufficiently large, so 
that the exact, stabilizing $n+m$ control kernels can be approximated sufficiently accurately (in 
specific sense) by the $\infty+\infty$ kernels constructed based on the continuum, fact which we 
prove. The proof relies on constructing proper sequences of backstepping kernels in $n$ and 
$m$, and showing that they converge to the continuum kernels, as $n, m \to \infty$. While the 
rationale of this design methodology stems from our earlier works \cite{HumBek25b, 
HumBek26}, establishing 
such an approximation property as $m\to \infty$ poses unique technical challenges, because the
number and form of the characteristics, along which the kernel equations are split into 
subdomains (where they are continuous), change with $m$. In the second, we 
construct controllers for stabilization of the large-scale hyperbolic system in the case where only 
macro measurements are available, i.e., in cases where only some average measurements 
originating from a macro (continuum) version of the original system are available, or when even 
only average measurements from that continuum system counterpart are available. To establish 
closed-loop 
stabilization we introduce a novel proof strategy in which we combine in a delicate manner the 
established solutions' and kernels' approximation property of the $n+m$ system by the 
respective 
$\infty+\infty$ continuum, with infinite-dimensional input-to-state stability (ISS) \cite{DasMir13, 
MirPri20} arguments. 

We furthermore provide a numerical example to illustrate stabilization of the continuum system 
itself, as well as to illustrate stabilization of the respective large-scale system, including 
verification of the limitations of our approach with respect to how large a large-scale system 
needs to be (i.e., how large $n$ and $m$ are required) for the controllers that employ continuum 
kernels to remain stabilizing. We also present a numerical example in which the objective is 
stabilization of a large-scale, $n+m$ hyperbolic system, when the control kernels are 
constructed as averaged continuum kernels and the available measurements are obtained as 
averaged continuum measurements.

\subsection{Organization}

The rest of the paper is organized as follows. In Section~\ref{sec:nm}, we derive a continuum 
approximation for large-scale $n+m$ systems  and formally show that the continuum 
$\infty+\infty$ system may 
approximate the $n+m$ system by establishing a connection between the respective systems' 
solutions (Theorem~\ref{thm:infsolappr}). In 
Section~\ref{sec:infbs}, we derive the (macro) backstepping control law for the class of (macro) 
$\infty+\infty$ 
hyperbolic systems and study stability of the closed-loop system constructing a Lyapunov 
functional (Theorem~\ref{thm:stab}); whereas the well-posedness of the respective continuum, 
backstepping kernel 
equations is 
established in Section~\ref{sec:infkwp} (Theorem~\ref{thm:infkwp}). In Section~\ref{sec:appr}, we 
develop micro-macro 
controllers for large-scale $n+m$ systems based on control kernels 
(Theorem~\ref{thm:infUappr1}) and/or measurements (Theorem~\ref{thm:infUappr2} and 
Proposition~\ref{thm:ave})
obtained on the basis of the $\infty+\infty$ continuum system. In 
Section~\ref{sec:numex}, we present numerical simulations to illustrate the theoretical results and 
the effectiveness of the presented control designs. Finally, Section~\ref{sec:conc} contains 
concluding remarks.

\subsection{Notation}

We use the standard notation $L^2(\Omega; \mathbb{R})$ for real-valued Lebesgue integrable 
functions on an arbitrary domain $\Omega$. Similarly, $L^\infty(\Omega; \mathbb{R}), C(\Omega; 
\mathbb{R}), C^1(\Omega; \mathbb{R})$ denote essentially bounded, continuous, and 
continuously differentiable functions, respectively, on $\Omega$. We occasionally use the 
shorthand $L^2$ when $\Omega$ is clear form the context. We introduce the continuum space 
$E_c = L^2([0,1]; L^2([0,1]; \mathbb{R}))$, equipped with the 
inner product
\begin{align}
\langle f_1, f_2 \rangle_{E_c} = \int\limits_0^1\int\limits_0^1 f_1(x,\zeta)f_2(x,\zeta)d\zeta dx.
\end{align}
Hence, $E_c^2$ can be viewed as the continuum limit of the 
space $E = L^2([0,1]; \mathbb{R}^{n+m})$ equipped with the inner product
\begin{align}
\left\langle \left(\begin{smallmatrix}
\mathbf{u}_1 \\ \mathbf{v}_1
\end{smallmatrix}\right), \left(\begin{smallmatrix}
\mathbf{u}_2 \\ \mathbf{v}_2
\end{smallmatrix}\right) \right\rangle_E & = \nonumber \\
\int\limits_0^1\frac{1}{n} \sum_{i=1}^n u_1^i(x)u_2^i(x)dx + \int\limits_0^1\frac{1}{m} 
\sum_{j=1}^m v_1^j(x)v_2^j(x)dx,
\end{align}
for some $n,m \in \mathbb{N}$, as $n,m \to \infty$. Moreover, we denote by $\mathcal{T}$ the 
triangular set
\begin{align}
	\mathcal{T} & = \left\{ (x,\xi) \in [0,1]^2: \xi \leq x\right\}.
\end{align}

For two normed spaces $Z, \mathcal{Z}$, we denote the space of bounded linear operators by 
$\mathcal{L}(Z,\mathcal{Z})$, and $\|\cdot\|_{\mathcal{L}(Z,\mathcal{Z})}$ denotes the 
corresponding operator norm. For $\mathcal{L}(Z,Z)$, we denote $\mathcal{L}(Z)$. Finally, we 
say that a system is exponentially stable on $Z$ if for any initial condition $z_0 \in Z$ the (weak) 
solution $z(t)$ of the system satisfies $\|z(t)\|_{Z} \leq Me^{-ct}\|z_0\|_{Z}$ for some $M,c > 0$ 
that are independent of $z_0$.

\section{Large-Scale Systems of $n+m$ PDEs and Convergence to an $\infty+\infty$ 
Continuum} \label{sec:nm}

\subsection{Large-Scale Systems of $n+m$ Hyperbolic PDEs}

The $n+m$ systems considered are of the form
\begin{subequations}
\label{eq:nmm}%
\begin{align}
	\mathbf{u}_t(t,x) + \pmb{\Lambda}(x)\mathbf{u}_x(t,x) 
	& = \frac{1}{n}\pmb{\Sigma}(x)\mathbf{u}(t,x) +
	\frac{1}{m}\mathbf{W}(x)\mathbf{v}(t,x), \\
	\mathbf{v}_t(t,x) - \mathbf{M}(x)\mathbf{v}_x(t,x)
	& = \frac{1}{n}\pmb{\Theta}(x)\mathbf{u}(t,x) + \frac{1}{m}\pmb{\Psi}(x)\mathbf{v}(t,x), 
\end{align}
\end{subequations}
with boundary conditions
\begin{align}
	\label{eq:nmmbc}%
	\mathbf{u}(t,0) & = \frac{1}{m}\mathbf{Q}\mathbf{v}(t,0), 
	& & \mathbf{v}(t,1) = \frac{1}{n}\mathbf{R}\mathbf{u}(t,1) + \mathbf{U}(t),
\end{align}
where $\mathbf{u} = (u_i)_{i=1}^n, \mathbf{v} = (v_j)_{j=1}^m$ are the states, $\mathbf{U} = 
(U_j)_{j=1}^m$ is the control input, and 
\begin{subequations}
\label{eq:nmmparam}
\begin{align}
\pmb{\Lambda} & = 
\operatorname{diag}(\lambda_1,\ldots,\lambda_n) \in C^1([0,1]; \mathbb{R}^{n\times n}), \\
 \mathbf{M} & =
\operatorname{diag}(\mu_1,\ldots,\mu_m) \in C^1([0,1]; \mathbb{R}^{m\times m}) \\
\pmb{\Sigma} & = (\sigma_{i,j})_{i,j=1}^n \in C([0,1]; \mathbb{R}^{n\times n}), \\
\mathbf{W} & = (w_{i,j})_{i=1,}^n{}_{j=1}^m \in C([0,1]; \mathbb{R}^{n\times m}), \\
\pmb{\Theta} & = (\theta_{j,i})_{j=1,}^m{}_{i=1}^n \in C([0,1]; \mathbb{R}^{m\times n}), \\
\pmb{\Psi} & = (\psi_{i,j})_{i,j=1}^m \in C([0,1]; \mathbb{R}^{m\times m}), \\
\mathbf{Q} & = (q_{i,j})_{i=1,}^n{}_{j=1}^m \in \mathbb{R}^{n\times m}, \\
 \mathbf{R} & = (r_{j,i})_{j=1}^m{}_{i=1}^n \in \mathbb{R}^{m\times n}.
\end{align}
\end{subequations}
As in \cite{AnfAamBook, 
HuLDiM16, HuLVaz19}, we make the following assumptions
on the parameters.
\begin{assumption}
\label{ass:nm}
The transport velocities in \eqref{eq:nmm} satisfy $\lambda_i(x) > 0$ for all $x \in 
[0,1]$ and $i = 1,\ldots,n$, and
\begin{equation}
\label{eq:munmass}
	\mu_1(x) > \mu_2(x) > \cdots > \mu_m(x) > 0,
\end{equation}
for all $x \in [0,1]$. Moreover, without loss of generality, we assume that $\psi_{j,j} = 0$ for all 
$j=1,\ldots,m$.\footnote{The diagonal terms of $\pmb{\Psi}$ can be canceled out through a 
change of variables (see, e.g., \cite[Sect. 3]{HuLVaz19}).}
\end{assumption}

\subsection{Continuum Approximation of Large-Scale $n+m$ Systems}

The goal of this subsection is to introduce a systematic approach for construction of an 
$\infty+\infty$ continuum system on $E_c^2$, 
which acts as a continuum approximation of the large-scale $n+m$ system \eqref{eq:nmm}, 
\eqref{eq:nmmbc}, in order to subsequently utilize it for control design for the $n+m$ system. We 
then specify the exact approximation properties as Theorem~\ref{thm:infsolappr}. As a first step 
towards this goal, we introduce a linear transform 
$\mathcal{F}_n$ (respectively, for $m$), which maps any vector $\mathbf{b}
=\left(b_i\right)_{i=1}^{n} \in \mathbb{R}^{n}$ into a step function in
$L^2\left([0,1];\mathbb{R}\right)$  as $\mathcal{F}_n\mathbf{b} = \sum_{i=1}^n b_i 
\chi_{((i-1)/n,i/n]}$, where $\chi_{((i-1)/n,i/n]}$ denotes the indicator function of the interval 
$((i-1)/n,i/n]$. Moreover, $\mathcal{F}_n$ is an isometry, i.e., it satisfies 
$\mathcal{F}_n^*\mathcal{F}_n = I_n$, where the adjoint $\mathcal{F}_n^{*}$ is given by
\begin{equation}
	\label{eq:Fns}
	\mathcal{F}_n^{*}h = \left( n\int\limits_{(i-1)/n}^{i/n}h(\zeta)d\zeta \right)_{i=1}^n,
\end{equation}
where each component is the mean value of any $h \in L^2([0,1];\mathbb{R})$ over the interval 
$[(i-1)/n,i/n]$.

We then apply the transform $\mathcal{F} = \operatorname{diag}(\mathcal{F}_n, \mathcal{F}_m)$ 
to \eqref{eq:nmm}, \eqref{eq:nmmbc} from the left to get
\begin{subequations}
\label{eq:nmmF}%
\begin{align}
	\mathcal{F}_n\mathbf{u}_t(t,x) + 
	\mathcal{F}_n\pmb{\Lambda}(x)\mathcal{F}_n^*\mathcal{F}_n\mathbf{u}_x(t,x) & = \nonumber \\
	\mathcal{F}_n\frac{1}{n}\pmb{\Sigma}(x)\mathcal{F}_n^*\mathcal{F}_n\mathbf{u}(t,x) +
	\mathcal{F}_n\frac{1}{m}\mathbf{W}(x)\mathcal{F}_m^*\mathcal{F}_m\mathbf{v}(t,x), \\
	\mathcal{F}_m\mathbf{v}_t(t,x) - 
	\mathcal{F}_m\mathbf{M}(x)\mathcal{F}_m^*\mathcal{F}_m\mathbf{v}_x(t,x)	& = \nonumber \\
	\mathcal{F}_m\frac{1}{n}\pmb{\Theta}(x)\mathcal{F}_n^*\mathcal{F}_n\mathbf{u}(t,x) + 	
	\mathcal{F}_m\frac{1}{m}\pmb{\Psi}(x)\mathcal{F}_m^*\mathcal{F}_m\mathbf{v}(t,x), 
	\end{align}
\end{subequations}
with boundary conditions
\begin{subequations}
\label{eq:nmmbcF}%
\begin{align}
		\mathcal{F}_n\mathbf{u}(t,0) & = 	
		\mathcal{F}_n\frac{1}{m}\mathbf{Q}\mathcal{F}_m^*\mathcal{F}_m\mathbf{v}(t,0), \\
		\mathcal{F}_m\mathbf{v}(t,1) & = 
		\mathcal{F}_m\frac{1}{n}\mathbf{R}\mathcal{F}_n^*\mathcal{F}_n\mathbf{u}(t,1) + 
		\mathcal{F}_m\mathbf{U}(t),
\end{align}
\end{subequations}
where we additionally use the isometry property of $\mathcal{F}_n$ and $\mathcal{F}_m$. Now, 
defining new state variables and input as
\begin{subequations}
\begin{align}
u^n(t,x,\cdot) & = 	\mathcal{F}_n\mathbf{u}(t,x), \qquad v^m(t,x,\cdot)  = 	
\mathcal{F}_m\mathbf{v}(t,x), \\
U^m(t,\cdot) & = \mathcal{F}_m\mathbf{U}(t),
\end{align}
\end{subequations}
the system \eqref{eq:nmmF}, \eqref{eq:nmmbcF} can be rewritten, for almost every $y,\eta \in 
[0,1]$, as
\begin{subequations}
\label{eq:infF}%
\begin{align}
	u_t^n(t,x,y) + \lambda^n(x,y)u_x^n(t,x,y)   & = \nonumber \\
	\resizebox{.92\columnwidth}{!}{$\displaystyle\int\limits_0^1\sigma^n(x,y,\zeta)u^n(t,x,\zeta)d\zeta
	+
	\int\limits_0^1W^{n,m}(x,y,\zeta)v^m(t,x,\zeta)d\zeta$}, \\
	v_t^m(t,x,\eta) - \mu^m(x,\eta)v^m_x(t,x,\eta)
	& =  \nonumber \\
	\resizebox{.92\columnwidth}{!}{$\displaystyle\int\limits_0^1\theta^{m,n}(x,\eta,\zeta)u^n(t,x,\zeta)d\zeta
	 + 
	\int\limits_0^1\psi^m(x,\eta,\zeta)v^m(t,x,\zeta)d\zeta$},
\end{align}
\end{subequations}
with boundary conditions 
\begin{subequations}
\label{eq:infbcF}%
\begin{align}
u^n(t,0,y) & = \int\limits_0^1Q^{n,m}(y,\zeta)v^m(t,0,\zeta)d\zeta, \\
v^m(t,1,\eta) & = \int\limits_0^1 R^{m,n}(\eta,\zeta)u^n(t,1,\zeta)d\zeta + U^m(t,\eta),
\end{align}
\end{subequations}
where, for all $x \in [0,1]$,
\begin{subequations}
\label{eq:infsparamF}
\begin{align}
& \mathcal{F} \begin{bmatrix}
\pmb{\Lambda}(x) & 0 \\ 0 & \mathbf{M}(x)
\end{bmatrix}  \mathcal{F}^*\begin{bmatrix}
u^n(t,x,\cdot) \\ v^m(t,x,\cdot)
\end{bmatrix} = \begin{bmatrix}
\lambda^n(x,\cdot)u^n(t,x,\cdot) \\ \mu^m(x,\cdot)v^m(t,x,\cdot)
\end{bmatrix}, \label{eq:infsparamF1} \\
& \mathcal{F} \begin{bmatrix}
\frac{1}{n}\pmb{\Sigma}(x) & \frac{1}{m}\mathbf{W}(x) \\ \frac{1}{n}\pmb{\Theta}(x) & 
\frac{1}{m}\pmb{\Psi}(x)
\end{bmatrix} \mathcal{F}^* \begin{bmatrix}
u^n(t,x,\cdot) \\ v^m(t,x,\cdot)
\end{bmatrix} = \nonumber \\
& \resizebox{!}{.053\textwidth}{$ \displaystyle \begin{bmatrix}
\int\limits_0^1 \sigma^n(x,\cdot,\zeta)u^n(t,x,\zeta)d\zeta & 
\int\limits_0^1 W^{n,m}(x,\cdot,\zeta)v^m(t,x,\zeta)d\zeta \\
\int\limits_0^1 \theta^{m,n}(x,\cdot,\zeta)u^n(t,x,\zeta)d\zeta & 
\int\limits_0^1 \psi^m(x,\cdot,\zeta)v^m(t,x,\zeta)d\zeta
\end{bmatrix}$}, \\
& \resizebox{!}{.051\textwidth}{$ \displaystyle\mathcal{F} \begin{bmatrix}
0 & \frac{1}{m}\mathbf{Q} \\ \frac{1}{n}\mathbf{R} & 0
\end{bmatrix} \mathcal{F}^* \begin{bmatrix}
u^n(t,1,\cdot) \\ v^m(t,0,\cdot)
\end{bmatrix} =
\begin{bmatrix}
\int\limits_0^1 Q^{n,m}(\cdot,\zeta)v^m(t,0,\zeta)d\zeta \\
\int\limits_0^1 R^{m,n}(\cdot,\zeta)u^n(t,1,\zeta)d\zeta
\end{bmatrix}$}.
\end{align}	
\end{subequations}

The system \eqref{eq:infF}, \eqref{eq:infbcF} is of the sought $\infty+\infty$ form, but considering 
that it is merely a representation of the original $n+m$ system \eqref{eq:nmm}, \eqref{eq:nmmbc} 
using step functions, one cannot argue yet that this is a continuum PDE system approximating 
\eqref{eq:nmm}, \eqref{eq:nmmbc}. However, based on \eqref{eq:infF}, \eqref{eq:infbcF}, we can 
now construct a continuum PDE system that acts as a continuum approximation of a large-scale 
$n+m$ system as follows
\begin{subequations}
\label{eq:inf}%
\begin{align}
u_t(t,x,y) + \lambda(x,y)u_x(t,x,y)   & = \nonumber \\
\int\limits_0^1\sigma(x,y,\zeta)u(t,x,\zeta)d\zeta +
\int\limits_0^1W(x,y,\zeta)v(t,x,\zeta)d\zeta, \label{eq:inf1} \\
v_t(t,x,\eta) - \mu(x,\eta)v_x(t,x,\eta)
& =  \nonumber \\
\int\limits_0^1\theta(x,\eta,\zeta)u(t,x,\zeta)d\zeta + 
\int\limits_0^1\psi(x,\eta,\zeta)v(t,x,\zeta)d\zeta, 
\label{eq:inf2}
\end{align}
\end{subequations}
with boundary conditions 
\begin{subequations}
\label{eq:infbc}%
\begin{align}
u(t,0,y) & = \int\limits_0^1Q(y,\zeta)v(t,0,\zeta)d\zeta, \\
v(t,1,\eta) & = \int\limits_0^1 R(\eta,\zeta)u(t,1,\zeta)d\zeta + U(t,\eta),
\end{align}
\end{subequations}
where the parameters are chosen such that, for a given $\varepsilon > 0$, they satisfy
\begin{subequations}
\label{eq:infappr}
\begin{align}
\max_{x\in[0,1]}\|\lambda^n(x,\cdot) - \lambda(x,\cdot)\|_{L^2([0,1]; \mathbb{R})} \qquad  
\nonumber \\
+ \max_{x\in[0,1]}\|\lambda_x^n(x,\cdot) - \lambda_x(x,\cdot)\|_{L^2([0,1]; \mathbb{R})} & < 
\varepsilon, \label{eq:infappr1}
\\
\max_{x\in[0,1]}\|\mu^m(x,\cdot) - \mu(x,\cdot)\|_{L^2([0,1]; \mathbb{R})} \qquad \nonumber \\ 
+\max_{x\in[0,1]}\|\mu_x^m(x,\cdot) - \mu_x(x,\cdot)\|_{L^2([0,1]; \mathbb{R})} & < \varepsilon,  
\\
\max_{x\in[0,1]}\|\sigma^n(x,\cdot) - \sigma(x,\cdot)\|_{L^2([0,1]^2; \mathbb{R})} & < 
\varepsilon,  
\\
\max_{x\in[0,1]}\|W^{n,m}(x,\cdot) - W(x,\cdot)\|_{L^2([0,1]^2; \mathbb{R})} & < \varepsilon,  \\
\max_{x\in[0,1]}\|\theta^{m,n}(x,\cdot) - \theta(x,\cdot)\|_{L^2([0,1]^2; \mathbb{R})} & < 
\varepsilon,  \\
\max_{x\in[0,1]}\|\psi^{m}(x,\cdot) - \psi(x,\cdot)\|_{L^2([0,1]^2; \mathbb{R})} & < \varepsilon,  \\
\|Q^{n,m}  - Q\|_{L^2([0,1]^2; \mathbb{R})} & < \varepsilon, \\
\|R^{m,n}  - R\|_{L^2([0,1]^2; \mathbb{R})}  & < \varepsilon.
\end{align}
\end{subequations}
In addition to the desired approximation accuracy as per 
\eqref{eq:infappr}, we make the following 
assumption about the parameters.\footnote{Assumption~\ref{ass:inf} contains the minimal 
assumptions about the parameters of \eqref{eq:inf}, \eqref{eq:infbc} for considering backstepping 
stabilization of this class of systems (see Section~\ref{sec:infbs}). Naturally, the parameters of 
\eqref{eq:inf}, \eqref{eq:infbc} 
can be also constructed such that they have more regularity, e.g., continuity in $\eta,\zeta$, but 
such 
additional regularity is not needed for studying backstepping control of \eqref{eq:inf},  
\eqref{eq:infbc}.}

\begin{assumption}
\label{ass:inf}
The parameters of \eqref{eq:inf}, \eqref{eq:infbc} are such that
$\lambda,\mu \in C^1([0,1]^2; \mathbb{R})$, $\sigma,W,\theta,\psi \in C([0,1];L^2([0,1]^2; 
\mathbb{R}))$, and $Q\in 
L^2([0,1]^2;\mathbb{R})$. Moreover, $\mu(x,\eta) 
> 0$ and $\lambda(x,y) > 0$ for all $x,y,\eta \in [0,1]$, and,
additionally,
\begin{equation}
\label{eq:muinfass}
\mu(x,\eta) > \mu(x,\zeta),
\end{equation}
for all $0 \leq \eta < \zeta \leq 1$ and $x \in [0,1]$. Finally, $\mu$ and $\psi$ are such that
\begin{equation}
\label{eq:mupsiass}
\max_{x\in [0,1]} \int\limits_0^1\int\limits_0^1 \left(\frac{\psi(x,\eta,\zeta)}{\mu(x,\eta) - 
\mu(x,\zeta)}\right)^2d\eta d\zeta < \infty.
\end{equation}
\end{assumption}

\begin{remark}
Conditions \eqref{eq:muinfass}, \eqref{eq:mupsiass}	are required for guaranteeing 
well-posedness of the resulting backstepping 
kernel equations, once we apply backstepping to the continuum system \eqref{eq:inf}, 
\eqref{eq:infbc} (see Sections~\ref{sec:infbs} and \ref{sec:infkwp}), where the assumption 
\eqref{eq:muinfass} about the $\mu$-velocities being 
ordered is consistent with the $n+m$ case \eqref{eq:munmass} (see, e.g., \cite[Sect. 
II]{HuLDiM16}). The assumption 
\eqref{eq:mupsiass}, on the other 
hand, is specific for the $\infty+\infty$ class of continuum systems, although it can be viewed as a 
counterpart of the $n+m$ assumption about the diagonal entries of $\psi$ being zero (cf. 
Assumption~\ref{ass:nm}), because in 
both cases such a condition stems from the boundary condition of the respective kernel 
equations. However, 
as $\psi \in C([0,1]; L^2([0,1]^2; \mathbb{R}))$, this assumption about the diagonal entries of 
$\psi$ does not translate as such to the continuum case, as the ``diagonal'' $\psi(x,\eta,\eta)$ 
may be 
ill-defined due to the line $\zeta = \eta$ being a measure zero subset of $(\eta,\zeta) \in 
[0,1]^2$. 
Hence, we have \eqref{eq:mupsiass} as a standing assumption.
\end{remark}

\begin{remark}
\label{rem:nmcont}
One option for obtaining functions $\lambda, \mu, \sigma, W, \theta, W, \psi, Q$, and $R$ that 
satisfy \eqref{eq:infappr}, based on the parameters \eqref{eq:nmmparam}, is to construct 
continuous (in all variables) functions, with the regularity of Assumption~\ref{ass:inf}, such that
\begin{subequations}
\label{eq:infcap}%
\begin{align}
\lambda(x,i/n) & = \lambda_i(x), \\
\mu(x, j/m) & = \mu_j(x), \\
\sigma(x,i/n,l/n) & = \sigma_{i,l}(x), \\
W(x,i/n,j/m) & = w_{i,j}(x), \\
\theta(x,j/m, i/n) & = \theta_{j,i}(x), \\
\psi(x,j/m, p/m) & = \psi_{j,p}(x), \label{eq:infcappsi} \\
Q(i/n,j/m) & = q_{i,j}, \\
R(j/m, i/n) & = r_{j,i},
\end{align}
\end{subequations}
for all $x \in [0,1]$, $i,l = 1,\ldots,n$, and $j,p = 1,\dots,m$, which can be done in infinitely many 
ways (see, e.g., \cite[Footnote 4]{HumBek25b}), but any such construction satisfies 
\eqref{eq:infappr} for any given $\varepsilon > 0$, when $n$ and $m$ are sufficiently large.
As \eqref{eq:infcap} do not characterize the continuum parameters uniquely, one also needs to 
ensure that the constructed continuum parameters additionally satisfy Assumption~\ref{ass:inf}, 
on the basis that the $n+m$ parameters satisfy Assumption~\ref{ass:nm}. This can be achieved, 
e.g., by taking $\mu$ of the form 
\begin{equation}
\mu(x,\eta) = \mu_m(x) + \sum_{\ell = 1}^{\widetilde{m}} a_\ell(x)(1-\eta)^\ell,
\end{equation}
for some $\widetilde{m} \in \mathbb{N}$, where $a_\ell(x) \geq 0$ for all $x \in [0,1]$ and $\ell \in 
\{1, 
\ldots,\widetilde{m}\}$ with $\displaystyle \min_{x \in [0,1]} \sum_{\ell = 1}^{\tilde{m}} a_\ell(x) > 
0$, which guarantees that $\mu(x,\eta) > 0$ and $\mu_y(x,\eta) < 0$ for all $x,\eta \in [0,1]$, so 
that 
\eqref{eq:muinfass} holds. Thereafter, $\psi$ can be taken of the form $\psi(x,\eta,\zeta) = 
\widetilde{\psi}(x,\eta,\zeta)(\mu(x,\eta) - \mu(x,\zeta))$, where $\widetilde{\psi}$ is 
constructed 
to satisfy
\begin{equation}
\widetilde{\psi}(x,j/m,p/m) = \frac{\psi_{j,p}(x)}{\mu_j(x) - \mu_p(x)}, \quad 1 \leq j \neq p \leq m,
\end{equation}
for all $x \in [0,1]$, so that \eqref{eq:mupsiass} and \eqref{eq:infcappsi} hold.
\end{remark}

\begin{remark}
If the parameters of the $n+m$ system are available or can be recast as expressions of $n$ and 
$m$, their continuum approximations can be taken as the limits of the respective 
sequences of step functions, e.g., (with reference to \eqref{eq:infsparamF1}) $\displaystyle 
\lambda = \lim_{n\to\infty}\lambda^n$, in which 
case $\varepsilon \to 0$ in \eqref{eq:infappr} as $n,m \to \infty$. In other words, one can derive, 
rather than construct, the parameters of the continuum system based on the parameters of the 
$n+m$ system. However, since the obtained 
continuum parameters need to satisfy Assumption~\ref{ass:inf}, additional conditions may be 
required to be imposed on $(\lambda_i)_{i=1}^n$ and $(\mu_j)_{j=1}^m$, to guarantee that their 
continuum limits are continuously differentiable functions in the ensemble variables (this is not an 
issue for the rest of the parameters whose limits are required to be only $L^2$ functions). For 
example, continuity of $\lambda$ can be guaranteed if $\lambda_i - \lambda_{i+1} \to 0$ for all $i 
= 1,\ldots,n-1$ as $n \to \infty$ and continuous differentiability if
$n(\lambda_{i+1} - 2\lambda_i + \lambda_{i-1}) \to 0$ for all $i = 2,\ldots,n-1$ as $n \to \infty$.
Note that when continuum approximations are constructed for given $n+m$ parameters along the 
lines of Remark~\ref{rem:nmcont}, Assumption~\ref{ass:inf} is always satisfied by construction.
\end{remark}

We end this section by establishing that \eqref{eq:inf}, \eqref{eq:infbc} is indeed a continuum 
approximation of \eqref{eq:nmmF}, \eqref{eq:nmmbcF} (and hence, of  \eqref{eq:nmm}, 
\eqref{eq:nmmbc}) in the sense that the solutions of the two systems remain arbitrarily close to 
each other on compact time intervals, provided that the respective parameters, initial conditions, 
and inputs are sufficiently close to 
each other.
\begin{theorem}
\label{thm:infsolappr}
Consider an $n+m$ system \eqref{eq:nmm}, \eqref{eq:nmmbc} with parameters $\pmb{\Lambda}, 
\mathbf{M}, \pmb{\Sigma}, \mathbf{W}, \pmb{\Theta}, \pmb{\Psi}, \mathbf{Q}, \mathbf{R}$ 
satisfying Assumption~\ref{ass:nm}, initial condition $(\mathbf{u}_0,\mathbf{v}_0) \in E$, and 
input $\mathbf{U} \in L^2([0, T]; \mathbb{R}^m)$ for an arbitrary, fixed $T > 0$. Construct a 
continuum system \eqref{eq:inf}, \eqref{eq:infbc} with parameters 
$\lambda,\mu,\sigma,W,\theta,\psi,Q,R$ satisfying Assumption~\ref{ass:inf} and 
\eqref{eq:infappr}, and equip  \eqref{eq:inf}, \eqref{eq:infbc} with initial conditions $u_0,v_0 \in 
E_c$ and input $U \in L^2([0,T]; L^2([0,1]; \mathbb{R}))$ such that
\begin{subequations}
\label{eq:initUappr}
\begin{align}
\left\|\mathcal{F}\left( \begin{smallmatrix}
\mathbf{u}_0 \\ \mathbf{v}_0
\end{smallmatrix}\right)- \left( \begin{smallmatrix}
u_0 \\ v_0 
\end{smallmatrix}\right) \right\|_{E_c^2} & < \varepsilon_{u,v}, \label{eq:initUappr1} \\
\|\mathcal{F}_m\mathbf{U} - U\|_{L^2([0,T]; L^2([0,1]; \mathbb{R}))} & < \varepsilon_U. 
\label{eq:initUappr2}
\end{align}
\end{subequations}
Then, there exists some $\delta_T  > 0$ depending continuously on $\varepsilon, 
\varepsilon_{u,v}$, 
and $\varepsilon_U$ such that
\begin{equation}
\label{eq:infsolappr}
\max_{t \in [0,T]} \left\| \mathcal{F}\left(\begin{smallmatrix}
\mathbf{u}(t) \\ \mathbf{v}(t)
\end{smallmatrix}\\\right) - \left(\begin{smallmatrix}
u(t) \\ v(t)
\end{smallmatrix}\right) \right\|_{E_c^2} < \delta_T,
\end{equation}
where $\delta_T\to 0$ as $\varepsilon, \varepsilon_{u,v}, \varepsilon_U \to 0$.

\begin{proof}
Firstly, the well-posedness of \eqref{eq:nmm}, \eqref{eq:nmmbc} under Assumption~\ref{ass:nm} 
has been established in \cite[Rem. 2]{HumBek26} based on \cite[Prop. A.1]{HumBek25b}, 
and the well-posedness of \eqref{eq:inf}, \eqref{eq:infbc} follows by the same arguments as for 
$\infty+1$ systems in \cite[Prop. B.1]{HumBek25b}. Hence, the (weak) solution to \eqref{eq:nmm}, 
\eqref{eq:nmmbc} satisfies $(\mathbf{u}, \mathbf{v}) \in C([0,T]; E)$ and the (weak) solution to 
\eqref{eq:inf}, \eqref{eq:infbc} satisfies $(u,v) \in C([0,T]; E_c^2)$. Consequently, the system  
\eqref{eq:infF}, \eqref{eq:infbcF} is well-posed and its (weak) solution is $(u^n(t), v^m(t)) := 
\mathcal{F}(\mathbf{u}(t), \mathbf{v}(t))$. In the following, we consider $(u^n,v^m)$ and 
\eqref{eq:infF}, \eqref{eq:infbcF} instead of $(\mathbf{u}, \mathbf{v})$ and \eqref{eq:nmm}, 
\eqref{eq:nmmbc}, as they are connected via the isometric transform $\mathcal{F}$.

Due to 
well-posedness of \eqref{eq:infF}, \eqref{eq:infbcF} and  \eqref{eq:inf},
\eqref{eq:infbc}, from \cite[Prop. 4.2.5]{TucWeiBook}, there exist families of linear operators 
$\mathbb{T}_t^{n,m}, \Phi_t^{n,m}$ and $\mathbb{T}_t, \Phi_t$, for $t \geq 0$, depending 
continuously on $\lambda^n, \mu^m, \sigma^n, W^{n,m}, \theta^{m,n}, \psi^m, Q^{n,m}, 
R^{m,n}$ and $\lambda, \mu, \sigma, W, \theta, \psi, Q, R$, respectively, such that the solutions to 
\eqref{eq:infF}, \eqref{eq:infbcF} and  \eqref{eq:inf}, \eqref{eq:infbc} can be written as
\begin{subequations}
\begin{align}
\left( \begin{smallmatrix}
u^n(t) \\ v^m(t)
\end{smallmatrix}\right) & = \mathbb{T}^{n,m}_t\left(\begin{smallmatrix}
u_0^n \\ v_0^m
\end{smallmatrix}\right) + \Phi_t^{n,m}U^m, \\
\left( \begin{smallmatrix}
u(t) \\ v(t)
\end{smallmatrix}\right)  & = \mathbb{T}_t\left(\begin{smallmatrix}
u_0 \\ v_0
\end{smallmatrix}\right) + \Phi_tU,
\end{align}
\end{subequations}
respectively. Computing the difference of the two solutions and using the triangle inequality gives, 
for each $t \in [0,T]$,
\begin{align}
\label{eq:solapprtpm}
\left\|\left( \begin{smallmatrix}
u^n(t) \\ v^m(t)
\end{smallmatrix}\right) - \left( \begin{smallmatrix}
u(t) \\ v(t)
\end{smallmatrix}\right) \right\|_{E_c^2} 
\leq \nonumber \\
\left\| (\mathbb{T}_t^{n,m} - \mathbb{T}_t
)\left(\begin{smallmatrix}
u_0^n \\ v_0^m
\end{smallmatrix}\right) \right\|_{E_c^2} 
+ \|\mathbb{T}_t\|_{\mathcal{L}(E_c^2)} \left\| \left( \begin{smallmatrix}
u_0^n \\ v_0^m
\end{smallmatrix}\right)  - \left( \begin{smallmatrix}
u_0 \\ v_0
\end{smallmatrix}\right) \right\|_{E_c^2} \nonumber \\
+ \|(\Phi_t^{n,m} - \Phi_t)U^m\|_{E_c^2} \nonumber \\
 + \|\Phi_T\|_{\mathcal{L}(L^2([0,T]; L^2([0,1]; 
\mathbb{R})),E_c^2)}\|U^m - U\|_{L^2([0,T]; L^2([0,1]; 
\mathbb{R}))},
\end{align}
where the first and third term become arbitrarily small when $\varepsilon$ in \eqref{eq:infappr} is 
sufficiently small, while the second and fourth term become arbitrarily small when 
$\varepsilon_{u,v}, \varepsilon_U$ in \eqref{eq:initUappr} are sufficiently small, as 
$\mathbb{T}_t$ and $\Phi_t$ are uniformly bounded on compact time intervals. Thus, 
\eqref{eq:infsolappr} follows after taking the maxima over $t \in [0,T]$ on both sides of 
\eqref{eq:solapprtpm}.
\end{proof}
\end{theorem}

\begin{remark}
\label{rem:nminf}
We note that 
$\varepsilon_{u,v}, \varepsilon_U$ can be made arbitrarily small, for example, by connecting the 
continuum initial 
conditions $(u_0,v_0)$ and input $U$ to ($\mathbf{u}_0, \mathbf{v}_0)$ and $\mathbf{U}$, 
respectively, analogously to \eqref{eq:infcap} (see 
also \cite[(28)]{HumBek25b}) and letting $n$ and $m$ be sufficiently large. Furthermore, 
$\varepsilon$ can be made arbitrarily small for sufficiently large $n$ and $m$, provided that the 
parameters of the continuum system are connected to the parameters of the $n+m$ system 
through \eqref{eq:infcap}.
\end{remark}

\section{Backstepping Stabilization of $\infty+\infty$ Systems: Macro Control with Macro 
Kernels and Measurements} \label{sec:infbs}

\subsection{Control Design via Backstepping}

The backstepping state feedback law to stabilize \eqref{eq:inf}, \eqref{eq:infbc} is of the form
\begin{align}
	\label{eq:infU}
	U(t,\eta) & = - \int\limits_0^1 R(\eta,\zeta)u(t,1,\zeta)d\zeta \nonumber \\
	& \qquad + \int\limits_0^1\int\limits_0^1 K(1,\xi,\eta,\zeta)u(t,\xi,\zeta)d\zeta d\xi \nonumber \\
	& \qquad + 	\int\limits_0^1\int\limits_0^1 L(1,\xi,\eta,\zeta)v(t,\xi,\zeta)d\zeta d\xi,
\end{align}
where $K,L \in L^\infty(\mathcal{T}; L^2([0,1]^2; \mathbb{R}))$ satisfy the kernel equations
\begin{subequations}
\label{eq:infk}
\begin{align}
\mu(x,\eta)K_x(x,\xi,\eta,\zeta) - K_\xi(x,\xi,\eta,\zeta)\lambda(\xi,\zeta) & \nonumber \\
- K(x,\xi,\eta,\zeta)\lambda_\xi(\xi,\zeta) 
& = \nonumber \\
\resizebox{.93\columnwidth}{!}{$\displaystyle\int\limits_0^1K(x,\xi,\eta,\chi)\sigma(\xi,\chi,\zeta)d\chi
 +
\int\limits_0^1L(x,\xi,\eta,\chi)\theta(\xi,\chi,\zeta)d\chi$}, \label{eq:infk1}
& \\
\mu(x,\eta)L_x(x,\xi,\eta,\zeta) + L_\xi(x,\xi,\eta,\zeta)\mu(\xi,\zeta) & \nonumber \\
+ L(x,\xi,\eta,\zeta)\mu_\xi(\xi,\zeta) 
& = \nonumber \\
\resizebox{.93\columnwidth}{!}{$\displaystyle 
\int\limits_0^1K(x,\xi,\eta,\chi)W(\xi,\chi,\zeta)d\chi 
+ 
\int\limits_0^1L(x,\xi,\eta,\chi)\psi(\xi,\chi,\zeta)d\chi$},  \label{eq:infk2}
&
\end{align}	
\end{subequations}
with boundary conditions, for almost all $\eta,\zeta \in [0,1]$,
\begin{subequations}
\label{eq:infkbc}
\begin{align}
L(x,x,\eta,\zeta) & = \frac{\psi(x,\eta,\zeta)}{ \mu(x,\zeta) - \mu(x,\eta)}, \label{eq:infkbc1} \\
K(x,x,\eta,\zeta) & = -\frac{\theta(x,\eta,\zeta)}{\lambda(x,\zeta) +  \mu(x,\eta)}, 
\label{eq:infkbc2}
\end{align}
for almost all $0 \leq \eta \leq \zeta \leq 1$,
\begin{align}
L(x,0,\eta,\zeta) & = 
\frac{1}{\mu(0,\zeta)}\int\limits_0^1K(x,0,\eta,\chi)\lambda(0,\chi)Q(\chi,\zeta)d\chi,\label{eq:infkbc3}
\end{align}	
and, for almost all $0 \leq \zeta < \eta \leq 1$,
\begin{align}
L(1,\xi,\eta,\zeta) & = l(\xi,\eta,\zeta),\label{eq:infkbc4} 
\end{align}	
\end{subequations}
where \eqref{eq:infkbc4} is an artificial boundary condition and
$l$ is chosen to be compatible with \eqref{eq:infkbc1} on $(x,\xi) = (1,1)$, which can be 
guaranteed, for example, by choosing
\begin{equation}
\label{eq:infl}
l(\xi,\eta,\zeta) = \frac{\psi(\xi,\eta,\zeta)}{\mu(\xi,\zeta) - \mu(\xi,\eta)}, \qquad \forall \zeta < 
\eta.
\end{equation}
We note that the kernel equations are understood in the sense that $K,L \in L^\infty(\mathcal{T}; 
L^2([0,1]^2; \mathbb{R}))$ (so that \eqref{eq:infkbc1} is legitimate under \eqref{eq:mupsiass}), as 
the boundary conditions on $(x,\xi) = (0,0)$ are (generally) 
over-determined (for $L$ on $\eta \leq \zeta$) because of \eqref{eq:infkbc1} and 
\eqref{eq:infkbc3}, 
\eqref{eq:infkbc2}, so that the equations cannot be interpreted pointwise in $(x,\xi) \in 
\mathcal{T}$. For more details, we refer to Section~\ref{sec:infkwp} and Appendix~\ref{app:ker}, 
where the derivation and well-posedness analysis of the kernel equations, respectively, are 
presented.

\subsection{Stability of the Closed-Loop System Under the Backstepping Control Law}

The stability result of the closed-loop system under the backstepping control law is presented in 
Theorem~\ref{thm:stab}. The proof is based on stability analysis of the target system resulting 
from the backstepping transformation, which essentially corresponds to the continuum limit of the 
respective $n+m$ target system as \cite[(12), (13)]{HuLDiM16}  $n,m \to \infty$ (see also 
Remark~\ref{rem:G}).

\begin{theorem} \label{thm:stab}
Under Assumption~ \ref{ass:inf}, the control law \eqref{eq:infU} exponentially stabilizes the 
system \eqref{eq:inf}, \eqref{eq:infbc} on $E_c^2$.

\begin{proof}
Firstly, the closed-loop system is well-posed, because we established the well-posedness of the 
open-loop system \eqref{eq:inf}, \eqref{eq:infbc} in Theorem~\ref{thm:infsolappr}, and hence, the 
well-posedness of the closed-loop system follows, e.g., by \cite[Cor. 5.5.1]{TucWeiBook}.

Secondly, given the backstepping kernels $K, L$ in the control law \eqref{eq:infU}, we introduce 
the following state transformation
\begin{subequations}
	\label{eq:infV}%
\begin{align}
\alpha(t,x,y) & = u(t,x,y), \label{eq:infV1} \\
\beta(t,x,\eta) & = v(t,x,\eta) - \int\limits_0^x
\int\limits_0^1K(x,\xi,\eta,\zeta)u(t,\xi,\zeta)d\zeta d\xi \nonumber  \\
& \qquad -\int\limits_0^x\int\limits_0^1 L(x,\xi,\eta,\zeta)v(t,\xi,\zeta)d\zeta d\xi, 
\label{eq:infV2}
\end{align}
\end{subequations}
such that the closed-loop system of \eqref{eq:inf}, \eqref{eq:infbc} with \eqref{eq:infU} is 
transformed to
\begin{subequations}
\label{eq:infts}%
\begin{align}
\alpha_t(t,x,y) + \lambda(x,y)\alpha_x(t,x,y) & = \nonumber \\
\int\limits_0^1 \sigma(x,y,\zeta)\alpha(t,x,\zeta)d\zeta +
\int\limits_0^1 W(x,y,\zeta)\beta(t,x,\zeta)d\zeta & \nonumber \\
+ \int\limits_0^1\int\limits_0^x C^+(x,\xi,y,\zeta)\alpha(t,\xi,\zeta)d\xi d\zeta & \nonumber \\
+ \int\limits_0^1\int\limits_0^x C^-(x,\xi,y,\zeta)\beta(t,\xi,\zeta)d\xi d\zeta, & \label{eq:infts1} \\
\beta_t(t,x,\eta) - \mu(x,\eta)\beta_x(t,x,\eta) & =  \nonumber \\
\int\limits_0^\eta G(x,\eta,\zeta)\beta(t,0,\zeta)d\zeta,
\label{eq:infts2}
\end{align}
\end{subequations}
where $C^+,C^- \in L^\infty(\mathcal{T}; L^2([0,1]^2; \mathbb{R})), G \in L^\infty([0,1]; 
L^2([0,1]^2;\mathbb{R}))$ with $G(x,\eta,\zeta) \equiv 0$ for $\zeta > \eta$ (so that the last 
integral of 
\eqref{eq:infts2} is taken only over $\zeta \in [0,\eta]$), with boundary conditions 
\begin{equation}
\label{eq:inftsbc}
\alpha(t,0,y) =  \int\limits_0^1Q(y,\zeta)\beta(t,0,\zeta)d\zeta, \qquad
\beta(t,1,\eta) \equiv 0.
\end{equation}
As the state transformation \eqref{eq:infV} is boundedly invertible by Lemma~\ref{lem:infkinv}, the 
well-posedness and stability of the target system \eqref{eq:infts}, \eqref{eq:inftsbc} are 
equivalent to those of the original closed-loop system \eqref{eq:inf}, \eqref{eq:infbc} with 
\eqref{eq:infU}.
	
Finally, for showing the exponential stability of the target system \eqref{eq:infts}, 
\eqref{eq:inftsbc} on $E_c^2$, consider a scalar $\delta > 0$ and a continuous function 
$D(\zeta) 
> 
0$ for all $\zeta \in [0,1]$, and construct a  candidate Lyapunov functional as
\begin{equation}
	\label{eq:lyap}
	\resizebox{.98\columnwidth}{!}{$\displaystyle V(t) = \int\limits_0^1 \int\limits_0^1 \left( 
	e^{-\delta
		x}\frac{\alpha^2(t,x,\zeta)}{\lambda(x,\zeta)} +
	e^{\delta x}\frac{D(\zeta)}{\mu(x,\zeta)}\beta^2(t,x,\zeta)\right)d\zeta dx$}.
\end{equation}
Computing $\dot{V}(t)$ and integrating by parts in $x$ gives 
\begin{align}
	\label{eq:lyapd}
	\dot{V}(t)
	& = \left[-e^{-\delta x}\|\alpha(t,x,\cdot)\|^2_{L^2} +
	e^{\delta x}\|\beta(t,x,\cdot)\|^2_{D}
	\right]_0^1 \nonumber \\
	& \quad - \delta\int\limits_0^1 \left(e^{-\delta
		x}\|\alpha(t,x,\cdot)\|_{L^2}^2 + e^{\delta
		x}\|\beta(t,x,\cdot)\|_{D}^2 \right)dx \nonumber \\
	& \quad \resizebox{.86\columnwidth}{!}{$ \displaystyle
		+ 2\int\limits_0^1 \int\limits_0^1 \int\limits_0^1e^{-\delta
			x}\frac{\alpha(t,x,\zeta)}{\lambda(x,\zeta)}\sigma(x,\zeta,\chi)\alpha(t,x,\chi)
		d\chi d\zeta dx$} \nonumber \\
	& \quad \resizebox{.86\columnwidth}{!}{$ \displaystyle+ 2\int\limits_0^1 \int\limits_0^1 
		\int\limits_0^1e^{-\delta
			x}\frac{\alpha(t,x,\zeta)}{\lambda(x,\zeta)}
		W(x,\zeta,\chi)\beta(t,x,\chi)d\chi d\zeta dx$} \nonumber \\
	& \quad \resizebox{.86\columnwidth}{!}{$ \displaystyle
		+ 2\int\limits_0^1 \int\limits_0^1 \int\limits_0^1
		\int\limits_0^x e^{-\delta x}
		\frac{\alpha(t,x,\zeta)}{\lambda(x,\zeta)}C^+(x,\xi,\zeta,\chi)\alpha(t,\xi,\chi)
		d\xi d\chi d\zeta dx$} \nonumber \\
	& \quad \resizebox{.86\columnwidth}{!}{$ \displaystyle+ 2\int\limits_0^1 
		\int\limits_0^1\int\limits_0^1 \int\limits_0^x
		e^{-\delta x}\frac{\alpha(t,x,\zeta)}{\lambda(x,\zeta)}
		C^-(x,\xi,\zeta,\chi)\beta(t,\xi,\chi)d\xi d\chi d\zeta dx$} \nonumber \\
	& \quad \resizebox{.86\columnwidth}{!}{$ \displaystyle
		+2 \int\limits_0^1\int\limits_0^1  \int\limits_0^\zeta e^{\delta 
			x}\frac{D(\zeta)}{\mu(x,\zeta)}\beta(t,x,\zeta)G(x,\zeta,\chi)\beta(t,0,\chi)d\chi d\zeta 
			dx$},
\end{align}
where $\|\cdot\|_{D}^2 = \left\langle \cdot, D\cdot
\right\rangle_{L^2}$ denotes the $D$-weighted inner product\footnote{We use the shorthand 
	notations $\|\alpha(t,x,\cdot)\|_{L^2}, \|\beta(t,x,\cdot)\|_{L^2}, \|\beta(t,x,\cdot)\|_D$ instead of 
	writing the integrals over $y$ explicitly. While this 
	is a slight abuse of notation (as these function may not be in $L^2$ for all $x$), these 
	expressions are valid appearing inside the integrals over $x$.}. Using the
following bounds
\begin{subequations}
	\label{eq:lyapb}%
	\begin{align}
		m_{\lambda} & = \min_{x,y \in [0,1]} \lambda(x,y), \label{eq:mlam} \\
		m_{\mu}  &= \min_{x,\eta\in
			[0,1]} \mu(x,\eta), \label{eq:mmu} \\
		M_{\sigma} & = \max_{x\in [0,1]} \left\|
		\int\limits_0^1\sigma(x,\cdot,\chi)d\chi
		\right\|_{L^2}, \label{eq:Ms} \\
		M_W & = \max_{x\in [0,1]} \left\|
		\int\limits_0^1 W(x,\cdot,\chi)d\chi
		\right\|_{L^2}, \label{eq:MW} \\
		M_{C^+} & = \esssup_{(x,\xi) \in \mathcal{T}} \left\|
		\int\limits_0^1C^+(x,\xi,\cdot,\chi)d\chi \right\|_{L^2},
		\\
		M_{C^-} & = \esssup_{(x,\xi) \in \mathcal{T}} \left\|
		\int\limits_0^1C^-(x,\xi,\cdot,\chi)d\chi \right\|_{L^2}, \\
		M_G & = \esssup_{x\in
			[0,1]}\left\| \int\limits_0^1 G(x,\cdot,\chi)d\chi \right\|_{L^2}, \\
		M_Q & = \left\|\int\limits_0^1 Q(\cdot,\chi)d\chi \right\|_{L^2}.
	\end{align}
\end{subequations}
the boundary conditions \eqref{eq:inftsbc}, the Cauchy-Schwarz
inequality, and $2\left\langle f,g \right\rangle_{L^2} \leq
\|f\|_{L^2}^2 +\|g\|^2_{L^2}$ for any $f,g \in L^2$, we can estimate \eqref{eq:lyapd}
as
\begin{align}
	\label{eq:lyapdest}
	\dot{V}(t)
	& \leq - \int\limits_0^1\left(D(\zeta) -
	M_Q^2\right)\beta^2(t,0,\zeta)d\zeta
	\nonumber \\
	& \quad - \delta\int\limits_0^1 \left(e^{-\delta
		x}\|\alpha(t,x,\cdot)\|_{L^2}^2 + e^{\delta
		x}\|\beta(t,x,\cdot)\|_{D}^2 \right)dx \nonumber \\
	& \quad + 2\int\limits_0^1e^{-\delta x}
	\frac{M_{\sigma}+M_{C^+}}{m_{\lambda}}\|\alpha(t,x,\cdot)\|_{L^2}^2dx
	\nonumber \\
	& \quad + \int\limits_0^1e^{-\delta x} \left(
	\frac{\|\alpha(t,x,\cdot)\|_{L^2}^2}{m_{\lambda}^2} +
	M_W^2\| \beta(t,x,\cdot)\|_{L^2}^2 \right)dx \nonumber \\
	& \quad + \int\limits_0^1e^{-\delta x} \left(
	\frac{\|\alpha(t,x,\cdot)\|_{L^2}^2}{m_{\lambda}^2} + M_{C^-}^2\|
	\beta(t,x,\cdot)\|_{L^2}^2  \right)dx \nonumber \\
	& \quad + M_G\int\limits_0^1 e^{\delta x}
	\frac{\|\beta(t,x,\cdot)\|_D^2}{m_\mu}dx \nonumber \\
	& \quad + \frac{M_Ge^{\delta}}{\delta m_{\mu}} \int\limits_0^1 \int\limits_\zeta^1 D(\chi)d\chi
	\beta^2(t,0,\zeta) d\zeta,
\end{align}
Now, $\dot{V}(t)$ can be guaranteed to be
negative definite by choosing $\delta$ and $D$, e.g., such that
\begin{subequations}
	\label{eq:lparam}%
	\begin{align}
		\delta & > \resizebox{.85\columnwidth}{!}{$\displaystyle \max \left\{
			\frac{2m_{\lambda}(M_{\sigma} + M_{C^+}) + 
				2}{m_{\lambda}^2}, \frac{M_W^2 + M_{C^-}^2}{m_{\mu}}+M_G\right\} $}, \\
		D(\zeta) & = c \exp\left(\frac{M_Ge^\delta}{\delta m_\mu}(1-\zeta)\right),
	\end{align}
\end{subequations}
for any $c > \max\{M_Q^2,1\}$, so that $D$ satisfies
\begin{equation}
	\label{eq:lparamD}
	D(\zeta) - \frac{M_Ge^{\delta}}{\delta m_{\mu}} \int\limits_\zeta^1 D(\chi)d\chi > M_Q^2,
\end{equation}
and $D(\zeta) \geq 1$ for all $\zeta \in [0,1]$, so that $\|f\|_D \geq \|f\|_{L^2}$ for any $f \in 
L^2([0,1]; 
\mathbb{R})$. Thus, the claim follows.
\end{proof}
\end{theorem}

\begin{remark}
\label{rem:G}	
The triangular structure of $G$ (in $(\eta,\zeta)$) is key in enabling the choice of the Lyapunov 
functional \eqref{eq:infV} for studying stability of the target system \eqref{eq:infts}, 
\eqref{eq:inftsbc}, particularly the weight $D$ satisfying \eqref{eq:lparamD}. Without 
the triangular structure of $G$, the condition \eqref{eq:lparamD} would become 
\begin{equation}
 D(\zeta) - \frac{M_Ge^{\delta}}{\delta m_{\mu}} \int\limits_0^1 D(\chi)d\chi > M_Q^2,
\end{equation}
which, by estimating $\displaystyle \int\limits_0^1 D(\chi)d\chi \geq \min_{\zeta \in [0,1]} 
D(\zeta)$, 
requires that
\begin{equation}
\label{eq:Dce}
\min_{\zeta \in [0,1]} D(\zeta)\left(1 - \frac{M_Ge^\delta}{\delta m_\mu}\right) > M_Q^2,
\end{equation}
which is impossible to hold if $\frac{M_Ge^\delta}{\delta m_\mu} > 1$, and hence, condition 
\eqref{eq:Dce} is dependent on the 
parameters (through $M_G$, where $G$ given in \eqref{eq:infGdef} depends on the 
kernels $K,L$, which, in turn, depend on the parameters) of \eqref{eq:inf}, \eqref{eq:infbc}.
\end{remark}

\section{Micro-Macro Control of Large-Scale Hyperbolic Systems} \label{sec:appr}

\subsection{Micro Control with Macro Kernels and Micro Measurements}

In this subsection we construct stabilizing control laws for the large-scale $n+m$ system applying 
$m$ control inputs, employing the $n+m$ measurements of the states of system \eqref{eq:nmm}, 
\eqref{eq:nmmbc}, and employing the continuum kernels \eqref{eq:infk}--\eqref{eq:infl}. 

Based on \eqref{eq:infU}, we can construct such stabilizing law as
\begin{align}
\label{eq:infUappr1}
\mathbf{U}(t) & = -\frac{1}{n}\mathbf{R}\mathbf{u}(t,1) +
 \frac{1}{n}\int\limits_0^1 \widetilde{\mathbf{K}}(1,\xi)\mathbf{u}(t,\xi)d\xi \nonumber   \\
& \qquad -\frac{1}{m}\int\limits_0^1 \widetilde{\mathbf{L}}(1,\xi)\mathbf{v}(t,\xi)d\xi,
\end{align}
where $\widetilde{\mathbf{K}} = (\widetilde{K}_{i,j})_{i=1}^m{}_{j=1}^n, \widetilde{\mathbf{L}} = 
(\widetilde{L}_{i,j})_{i,j=1}^m$ are obtained through mean-value sampling of the continuum 
kernels as
\begin{subequations}
\label{eq:infUappr1k}
\begin{align}
\widetilde{K}_{i,j}(1,\xi) & = 
nm\int\limits_{(i-1)/m}^{i/m}\int\limits_{(j-1)/n}^{j/n}K(1,\xi,\eta,\zeta)d\zeta 
d\eta, \label{eq:infUappr1ka}
\\
\widetilde{L}_{i,j}(1,\xi) & = 
m^2\int\limits_{(i-1)/m}^{i/m}\int\limits_{(j-1)/m}^{j/m}L(1,\xi,\eta,\zeta)d\zeta 
d\eta,
\end{align}
\end{subequations}
for almost all $\xi \in [0,1]$. When $n,m$ are sufficiently large, the control law 
\eqref{eq:infUappr1}, \eqref{eq:infUappr1k} exponentially stabilizes the closed-loop system, as 
formally stated in the following theorem.

\begin{theorem}
\label{thm:infUappr1}
Under Assumption~\ref{ass:nm}, and provided that the continuum parameters are constructed 
such that Assumption~\ref{ass:inf} holds and \eqref{eq:infcap} is satisfied, the 
continuum-based control law \eqref{eq:infUappr1}, \eqref{eq:infUappr1k} exponentially stabilizes 
the $n+m$ system \eqref{eq:nmm}, \eqref{eq:nmmbc}, provided that $n,m$ are sufficiently large.

\begin{proof}
We begin with the exact $n+m$ kernel equations for $\mathbf{K} =
(K_{i,j})_{i=1}^m{}_{j=1}^n, \mathbf{L} = (L_{i,j})_{i,j=1}^m$, given by
\begin{subequations}
 \label{eq:nmk}
\begin{align}
\mathbf{M}(x)\mathbf{K}_x(x,\xi) - \mathbf{K}_{\xi}(x,\xi)\pmb{\Lambda}(\xi) -
  \mathbf{K}(x,\xi)\pmb{\Lambda}'(\xi)
  & = \nonumber \\  
 \mathbf{K}(x,\xi)\pmb{\Sigma}(\xi) +
\mathbf{L}(x,\xi)\pmb{\Theta}(\xi), \label{eq:nmk1} \\
  \mathbf{M}(x)\mathbf{L}_x(x,\xi) + \mathbf{L}_{\xi}(x,\xi)\mathbf{M}(\xi) +
  \mathbf{L}(x,\xi)\mathbf{M}'(\xi)
  & = \nonumber \\
  \mathbf{K}(x,\xi)\mathbf{W}(\xi) + \mathbf{L}(x,\xi)\pmb{\Psi}(\xi), \label{eq:nmk2}
\end{align}
\end{subequations}
with boundary conditions
\begin{subequations}
  \label{eq:nmkbc}%
  \begin{align}
0 & = \mathbf{K}(x,x)\pmb{\Lambda}(x) +
                   \mathbf{M}(x)\mathbf{K}(x,x) + \pmb{\Theta}(x), \label{eq:nmkbc1} \\
0 & = \mathbf{M}(x)\mathbf{L}(x,x) - \mathbf{L}(x,x)\mathbf{M}(x)+\pmb{\Psi}(x), 
\label{eq:nmkbc2} \\
  L_{i,j}(x,0)
  & = \frac{1}{\mu_j(0)} \sum_{\ell=1}^n
                 \lambda_{\ell}(0)K_{i,\ell}(x,0)Q_{\ell,j}, \quad \forall i
    \leq j,  \label{eq:nmkbc3} \\
L_{i,j}(1,\xi) & = l_{i,j}(\xi), \quad \forall j < i,  \label{eq:nmkbc4}
\end{align}
\end{subequations}
where $l_{i,j}$ are chosen such that the $\mathbf{L}$ kernels satisfy
a compatibility condition on $(x,\xi) = (1,1)$, e.g.,
\begin{equation}
  \label{eq:nml}
  l_{i,j}(\xi) = - \frac{\psi_{i,j}(\xi)}{\mu_i(\xi) - \mu_{j}(\xi)},
  \quad \quad \forall j < i,
\end{equation}
is a viable choice. In order to compare the $n+m$ kernels to the
$\infty+\infty$ kernels, we apply the transform $\mathcal{F}_m$ to
\eqref{eq:nmk}, \eqref{eq:nmkbc} from the left and
$\mathcal{F}_n^{*},\mathcal{F}_m^{*}$ from the right to obtain
\begin{subequations}
 \label{eq:nmkF}
\begin{align}
 \resizebox{.95\columnwidth}{!}{$\displaystyle 
 \mathcal{F}_m\mathbf{M}(x)\mathcal{F}_m^{*}\mathcal{F}_m\mathbf{K}_x(x,\xi)\mathcal{F}_n^{*}
  -\mathcal{F}_{m} \mathbf{K}_{\xi}(x,\xi) \mathcal{F}_n^{*}
  \mathcal{F}_n\pmb{\Lambda}(\xi)\mathcal{F}_n^{*}$}
  & \nonumber \\
  - \mathcal{F}_m\mathbf{K}(x,\xi) \mathcal{F}_n^{*} 
  \mathcal{F}_n\pmb{\Lambda}'(\xi)\mathcal{F}_n^{*}
  & = \nonumber \\  
 \mathcal{F}_m\mathbf{K}(x,\xi)\mathcal{F}_n^{*}\mathcal{F}_n\pmb{\Sigma}(\xi)\mathcal{F}_n^{*}
  +
\mathcal{F}_m\mathbf{L}(x,\xi)\mathcal{F}_m^{*}\mathcal{F}_m\pmb{\Theta}(\xi)\mathcal{F}_n^{*},
  \\
 \resizebox{.95\columnwidth}{!}{$\displaystyle\mathcal{F}_m\mathbf{M}(x)\mathcal{F}_{m}^{*}\mathcal{F}_m\mathbf{L}_x(x,\xi)\mathcal{F}_{m}^{*}
  + \mathcal{F}_m\mathbf{L}_{\xi}(x,\xi)\mathcal{F}_{m}^{*}\mathcal{F}_m\mathbf{M}(\xi) 
  \mathcal{F}_{m}^{*}$}
  & \nonumber \\ 
  +\mathcal{F}_m\mathbf{L}(x,\xi)\mathcal{F}_{m}^{*}\mathcal{F}_m\mathbf{M}'(\xi)\mathcal{F}_{m}^{*}
  & = \nonumber \\
  \mathcal{F}_m\mathbf{K}(x,\xi)\mathcal{F}_n^{*} 
  \mathcal{F}_n\mathbf{W}(\xi)\mathcal{F}_{m}^{*} +
  \mathcal{F}_m\mathbf{L}(x,\xi)\mathcal{F}_{m}^{*}\mathcal{F}_m\pmb{\Psi}(\xi)\mathcal{F}_{m}^{*},
\end{align}
\end{subequations}
with boundary conditions
\begin{subequations}
  \label{eq:nmkbcF}%
  \begin{align}
 \resizebox{.93\columnwidth}{!}{$\displaystyle\mathcal{F}_m\mathbf{K}(x,x)\mathcal{F}_n^{*} 
 \mathcal{F}_n\pmb{\Lambda}(x)\mathcal{F}_n^{*} +
     \mathcal{F}_m\mathbf{M}(x)\mathcal{F}_m^{*}\mathcal{F}_m\mathbf{K}(x,x)\mathcal{F}_n^{*}$}
    & = \nonumber \\
     \mathcal{F}_m\pmb{\Theta}(x)\mathcal{F}_n^{*}, \\
\resizebox{.93\columnwidth}{!}{$\displaystyle 
\mathcal{F}_m\mathbf{M}(x)\mathcal{F}_m^{*}\mathcal{F}_m\mathbf{L}(x,x)\mathcal{F}_m^{*} -
    \mathcal{F}_m\mathbf{L}(x,x)\mathcal{F}_m^{*}\mathcal{F}_m\mathbf{M}(x))\mathcal{F}_m^{*}$}
    & = \nonumber \\
    \mathcal{F}_m\pmb{\Psi}(x)\mathcal{F}_m^{*},
\end{align}
\end{subequations}
which are of the form of the respective $\infty+\infty$ kernel equations for
piecewise constant parameters defined in \eqref{eq:infsparamF}. Respectively, the boundary 
conditions \eqref{eq:nmkbc3}, \eqref{eq:nmkbc4} get transformed into piecewise boundary 
conditions in $(\eta,y)$ as
\begin{subequations}
\label{eq:nmkbcF2}
\begin{align}
\label{eq:nmkbcF21}
 \quad L^m(x,0,\eta,\zeta) & = \nonumber \\
\frac{1}{\mu^m(0,\zeta)}\int\limits_0^1K^{m,n}(x,0,\eta,\chi)\lambda^n(0,\chi)Q^{n,m}(\chi,\zeta)d\chi,
\end{align}
for all $(\eta,\zeta) \in ((i-1)/m,i/m]\times(((j-1)/m,j/m]$ with $1 \leq i\leq j \leq m$, and 
\begin{align}
\label{eq:nmkbcF22}
L^m(1,\xi,\eta,\zeta) = l^{m}(\xi,\eta,\zeta),
\end{align}
\end{subequations}
for all $(\eta,\zeta) \in ((i-1)/m,i/m]\times(((j-1)/m,j/m]$ with $1 \leq j < i \leq m$, where we use the 
notation\footnote{Note that $\mathcal{F}_m\mathbf{K}\mathcal{F}_n^*$ (resp. 
$\mathcal{F}_m\mathbf{L}\mathcal{F}_m^*$) is an integral operator, i.e.,
$\mathcal{F}_m\mathbf{K}\mathcal{F}_n^*h = \int\limits_0^1 
K^{m,n}(x,\xi,\eta,\zeta)h(\zeta)d\zeta$, for any $h \in L^2([0,1]; \mathbb{R})$.}
\begin{subequations}
\begin{align}
K^{m,n}(x,\xi,\eta,\zeta) & = K_{i,j}(x,\xi), & &  \eta \in ((i-1)/m,i/m],\nonumber \\
& & & \quad  \zeta \in ((j-1)/n,j/n], \\
L^{m}(x,\xi,\eta,\zeta) & = L_{i,j}(x,\xi), & & \eta \in ((i-1)/m,i/m], \nonumber \\
& & & \quad \zeta \in ((j-1)/m,j/m]
\end{align}
\end{subequations}
for almost every $(x,\xi) \in \mathcal{T}$. The $(\eta,\zeta)$-domains of the boundary 
conditions 
\eqref{eq:nmkbcF2} are illustrated in Fig.~\ref{fig:bct}, where one can see that in the limit case 
$m\to\infty$ the respective $(\eta,\zeta)$-domains of the continuum boundary conditions 
\eqref{eq:infkbc3}, \eqref{eq:infkbc4} are recovered. Also formally, the domains 
of \eqref{eq:nmkbcF21} and \eqref{eq:infkbc3} (respectively, \eqref{eq:nmkbcF22} and 
\eqref{eq:infkbc4}) differ in $(\eta,\zeta)$ by a measure of $\frac{m}{2}\frac{1}{m^2} = 
\frac{1}{2m}$, which vanishes as $m \to\infty$, so that the $\infty+\infty$ kernel equations 
\eqref{eq:infk}, \eqref{eq:infkbc} are recovered when $n,m\to\infty$, provided that
\eqref{eq:infcap} holds.

\begin{figure}[!htb]
\begin{center}
\includegraphics[scale=.69]{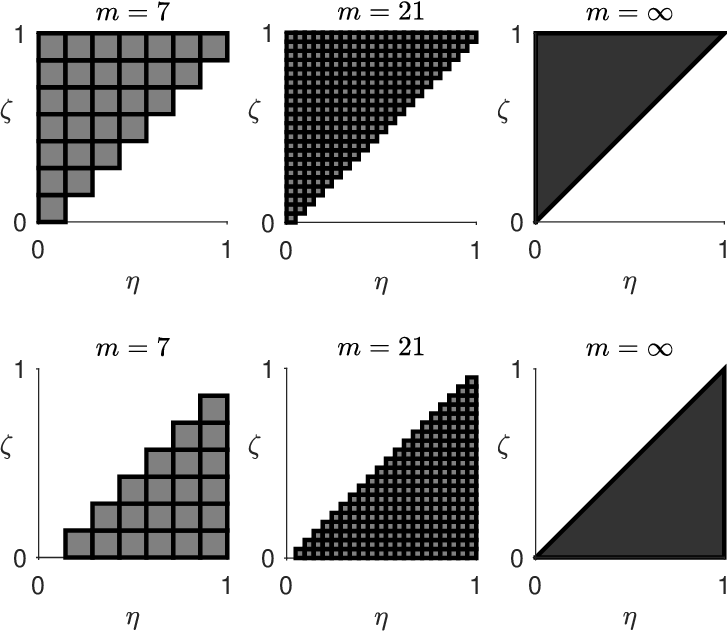}
\end{center}
\caption{Illustrations for the domains of the boundary conditions \eqref{eq:nmkbcF21} (upper row) 
and \eqref{eq:nmkbcF22} (lower row) for $m=7,21$, and the limits as $m\to\infty$, corresponding 
to the domains of \eqref{eq:infkbc3}, \eqref{eq:infkbc4} in $(\eta,\zeta)$.}
\label{fig:bct}
\end{figure}

The kernel equations \eqref{eq:infk}, \eqref{eq:infkbc} and \eqref{eq:nmk}, 
\eqref{eq:nmkbc} are well-posed by Theorem~\ref{thm:infkwp} and \cite[Sect. VI]{HuLDiM16}, 
respectively, and the solution $(\mathbf{K}, \mathbf{L})$ to \eqref{eq:nmk}, \eqref{eq:nmkbc} 
satisfies \eqref{eq:nmkF}, \eqref{eq:nmkbcF} by construction. Due to well-posedness, the 
solutions to the kernel equations depend continuously on the respective parameters, and hence, 
as $\varepsilon$ in \eqref{eq:infappr} becomes sufficiently small, the solutions to the $n+m$ and 
$\infty+\infty$ kernel equations satisfy
\begin{subequations}
\label{eq:infkcomp}
\begin{align}
\esssup_{(x,\xi) \in \mathcal{T}} \| K^{m,n}(x,\xi,\cdot,\cdot) - 
K(x,\xi,\cdot,\cdot)\|_{L^2([0,1]^2; \mathbb{R})} & \leq \delta_\varepsilon, \\
\esssup_{(x,\xi) \in \mathcal{T}} \| L^m(x,\xi,\cdot,\cdot) - 
L(x,\xi,\cdot,\cdot)\|_{L^2([0,1]^2; \mathbb{R})} & \leq \delta_\varepsilon,
\end{align}
\end{subequations}
where $\delta_\varepsilon$ depends continuously on $\varepsilon$ such that $\delta_\varepsilon 
\to 0$ as $\varepsilon \to 0$. Thus, when $n,m$ are sufficiently large,  \eqref{eq:infcap} implies 
that $\varepsilon$ in \eqref{eq:infappr} becomes arbitrarily small, and hence, $\delta_\varepsilon$ 
in \eqref{eq:infkcomp} becomes arbitrarily small as well.

Now, comparing the approximate control law \eqref{eq:infUappr1} with the exact $n+m$ law given 
by
\begin{align}
\label{eq:nmU}
\mathbf{U}_{\mathrm{e}}(t) & = -\frac{1}{n}\mathbf{R}\mathbf{u}(t,1) +
 \frac{1}{n}\int\limits_0^1 \mathbf{K}(1,\xi)\mathbf{u}(t,\xi)d\xi \nonumber   \\
& \qquad +\frac{1}{m}\int\limits_0^1 \mathbf{L}(1,\xi)\mathbf{v}(t,\xi)d\xi,
\end{align}
the first term is the same, so that it remains to estimate the kernels approximation errors. We get 
for the $\mathbf{K}, \widetilde{\mathbf{K}}$ terms (respectively for $\mathbf{L}, 
\widetilde{\mathbf{L}}$)
\begin{align}
\frac{1}{n}\int\limits_0^1 (\mathbf{K}(1,\xi) - 
\widetilde{\mathbf{K}}(1,\xi))\mathbf{u}(t,\xi)d\xi & = \nonumber \\
\int\limits_0^1 \mathcal{F}_m^*\left(\mathcal{F}_m\mathbf{K}(1,\xi) \mathcal{F}_n^* - 
\mathcal{F}_m\widetilde{\mathbf{K}}(1,\xi)\mathcal{F}_n^*\right)\mathcal{F}_n\mathbf{u}(t,\xi)d\xi,
\end{align}
where, for any $h \in L^2([0,1]; \mathbb{R})$,  
$\mathcal{F}_m\widetilde{\mathbf{K}}(1,\xi)\mathcal{F}_n^*h$ approximates 
$\int\limits_0^1K(1,\xi,\cdot,\zeta)h(\zeta)d\zeta$ for almost every $\xi \in [0,1]$ through the 
mean-value 
approximation \eqref{eq:infUappr1ka}. Hence,
using the triangle inequality and $\|\mathcal{F}^*\|_{\mathcal{L}(E_c^2, E)} = 1$, we get
\begin{align}
\left\|\frac{1}{n}\int\limits_0^1 (\mathbf{K}(1,\xi) - 
\widetilde{\mathbf{K}}(1,\xi))\mathbf{u}(t,\xi)d\xi\right\|_{\mathbb{R}^m} & \leq \nonumber \\
 \nonumber \\
\left(\resizebox{.9\columnwidth}{!}{$\displaystyle \sup_{h \in L^2:  \|h\| 
=1}\left\|\int\limits_0^1 \left( 
\int\limits_0^1K(1,\xi,\cdot,\zeta)h(\zeta)d\zeta - 
\mathcal{F}_m\widetilde{\mathbf{K}}(1,\xi)\mathcal{F}_n^*h\right)d\xi\right\|_{L^2([0,1]; 
\mathbb{R})}$}
\right. \nonumber \\
\left.+ \delta_\varepsilon\right)\frac{1}{\sqrt{n}}\|\mathbf{u}(t)\|_{\mathbb{R}^n},
\end{align}
where the mean-value approximation error becomes arbitrarily small when $n,m$ are sufficiently 
large due to step functions being dense in $L^2$ (see, e.g., \cite[Sect. 1.3.5]{TaoBook11}), and 
$\delta_\varepsilon$ from \eqref{eq:infkcomp} becomes arbitrarily small as $n,m$ are sufficiently 
large by the preceding arguments. Thus, the effect of the kernels approximation error tends to 
zero as $n,m$ tend to 
infinity, so that the stated exponential stability result follows by 
combining the exponential stability of the closed-loop under the exact backstepping control law 
(see, e.g., \cite[Thm 3.4]{HuLDiM16}) and robustness of exponential stability under sufficiently 
small, admissible perturbations (see, e.g., \cite[Prop. A.2]{HumBek25b}), when $n,m$ are 
sufficiently large.
\end{proof}
\end{theorem}

\subsection{Micro Control with Macro Kernels and Macro Measurements}

In this subsection we still consider that $m$ different controls are applied in \eqref{eq:nmmbc} 
and that the macro kernels, i.e., the kernels constructed based on the continuum system in 
Section~\ref{sec:infbs}, are employed. The difference is that we employ continuum-based 
measurements instead of the $n+m$ (micro) measurements (exactly) corresponding to the states 
of the $n+m$ system.

In order to present the respective control law, we introduce macro measurements $(\tilde{u}, 
\tilde{v})$ that approximate the full state information of the $n+m$ system as follows
\begin{subequations}
\label{eq:uveps}
\begin{align}
\sup_{t \geq 0}\left\| \left(\begin{smallmatrix}
\tilde{u}(t) \\ \tilde{v}(t)
\end{smallmatrix}\right) - \mathcal{F} \left(\begin{smallmatrix}
\mathbf{u}(t) \\ \mathbf{v}(t)
\end{smallmatrix}\right) \right\|_{E_c^2} & \leq \tilde{\varepsilon}_1,  \label{eq:uveps1} \\
\esssup_{t \geq 0}\|\tilde{u}(t,1,\cdot) - \mathcal{F}_n \mathbf{u}(t,1)\|_{L^2([0,1]; \mathbb{R})} & 
\leq \tilde{\varepsilon}_2, \label{eq:uveps2}
\end{align}
\end{subequations}
where $\tilde{\varepsilon}_1, \tilde{\varepsilon}_2 > 0$ determine the approximation accuracy of 
the macro measurement. The control law is then defined as
\begin{align}
\label{eq:infUappr2}
\mathbf{U}(t) & = -\int\limits_{0}^1\widetilde{\mathbf{R}}(\zeta)\tilde{u}(t,1,\zeta)d\zeta \nonumber 
\\
& \qquad + \int\limits_0^1\int\limits_0^1 \mathbf{K}(1,\xi,\zeta)\tilde{u}(t,\xi,\zeta)d\xi d\zeta 
\nonumber \\
& \qquad + \int\limits_0^1 \int\limits_0^1 \mathbf{L}(1,\xi,\zeta)\tilde{v}(t,\xi,\zeta)d\xi d\zeta,
\end{align}
with $\widetilde{\mathbf{R}} = (\widetilde{R}_i)_{i=1}^m$ given by
\begin{equation}
\label{eq:Rtilde}
\widetilde{R_i}(\zeta) = m\int\limits_{(i-i)/m}^{i/m}R(\chi,\zeta)d\chi,
\end{equation}
and $\mathbf{K} = 
(K_i)_{i=1}^m, 
\mathbf{L} = (L_i)_{i=1}^m$ given by
\begin{subequations}
\label{eq:infUappr2k}
\begin{align}
K_{i}(1,\xi,\zeta) & = m\int\limits_{(i-1)/m}^{i/m} K(1,\xi,\chi,\zeta)d\chi, 
\\
L_i(1,\xi,\zeta) & = m\int\limits_{(i-1)/m}^{i/m} L(1,\xi,\chi,\zeta)d\chi,
\end{align}
\end{subequations}
for almost all $\xi,\zeta \in [0,1]$. 

\begin{theorem}
\label{thm:infUappr2}
Let Assumption~\ref{ass:nm} hold and consider continuum parameters constructed 
such that Assumption~\ref{ass:inf} holds and \eqref{eq:infcap} is satisfied. Moreover, assume that
macro measurements $(\tilde{u}, \tilde{v})$ are available such that \eqref{eq:uveps} is satisfied for 
some $\tilde{\varepsilon}_1, \tilde{\varepsilon}_2 > 0$. Then, provided that $m,n$ are sufficiently 
large, there exist some $\tilde{M}, \tilde{\omega}, \tilde{H} > 0$ such 
that  the solution to $n+m$ system \eqref{eq:nmm}, \eqref{eq:nmmbc} under the control law 
\eqref{eq:infUappr2}--\eqref{eq:infUappr2k} satisfies
\begin{align}
\label{eq:macrostab}
\left\|\left(\begin{smallmatrix}
\mathbf{u}(t) \\ \mathbf{v}(t) \end{smallmatrix}\right)\right\|_{E} & \leq 
\tilde{M}e^{-\tilde{\omega}t}\left\|\left(\begin{smallmatrix}
\mathbf{u}_0 \\ \mathbf{v}_0 \end{smallmatrix}\right)\right\|_{E} \nonumber \\
& \quad + 
\tilde{H}\left(\tilde{M}_R\tilde{\varepsilon}_2 + (\tilde{M}_K + \tilde{M}_L)\tilde{\varepsilon}_1 
\right),
\end{align}
for any initial conditions $\left(\begin{smallmatrix} \mathbf{u}_0 \\ \mathbf{v}_0 
\end{smallmatrix}\right) \in E$, where we denote
\begin{subequations}
\begin{align}
\tilde{M}_R & = \|\mathbf{R}\|_{L^2([0,1]; \mathbb{R}^m)}, \\
\tilde{M}_K & = \left\|\int\limits_0^1 \mathbf{K}(1,\xi,\cdot)d\xi\right\|_{L^2([0,1]; \mathbb{R}^m)}, \\
\tilde{M}_L & = \left\|\int\limits_0^1 \mathbf{L}(1,\xi,\cdot)d\xi\right\|_{L^2([0,1]; \mathbb{R}^m)}. 
\end{align}
\end{subequations}

\begin{proof}
The basis of the proof is splitting the terms in the control law \eqref{eq:infUappr2} as (and 
analogously for $\mathbf{K}, \mathbf{L}$)
\begin{align}
\label{eq:apprtmp}
\int\limits_0^1\widetilde{\mathbf{R}}(\zeta)\tilde{u}(t,1,\zeta)d\zeta & = 
\frac{1}{n}\mathbf{R}\mathbf{u}(t,1)  +
\frac{1}{n}\left(\widetilde{\mathbf{R}} - \mathbf{R}\right)\mathbf{u}(t,1) \nonumber \\
& \quad + \int\limits_0^1 \widetilde{\mathbf{R}}(\zeta)\left( 
\tilde{u}(t,1,\zeta) - u^n(t,1,\zeta) \right)d\zeta,
\end{align}
where we denote $\widetilde{\mathbf{R}} = (\tilde{r}_{j,i})_{j=1}^m{}_{i=1}^n$ with 
$\tilde{r}_{j,i} = \int\limits_{(i-1)/n}^{i/n}\widetilde{R}_j(\zeta)d\zeta$, so that 
$\int\limits_0^1\widetilde{\mathbf{R}}(\zeta)u^n(t,1,\zeta)d\zeta = 
\frac{1}{n}\widetilde{\mathbf{R}}\mathbf{u}(t,1)$.
The first term of \eqref{eq:apprtmp} coincides with the respective term in the exact stabilizing 
control law \eqref{eq:nmU}, the second term is related to continuum approximation errors of the 
parameters/kernels, and the last term is related to macro measurement errors.

By analogous arguments to the proof of Theorem~\ref{thm:infUappr1}, the continuum 
approximation errors of $\widetilde{\mathbf{R}}, \mathbf{K}$, and $\mathbf{L}$ become arbitrarily 
small as $n,m$ are sufficiently large. Ignoring the macro measurement errors momentarily, 
exponential stability under the exact control law would be preserved despite sufficiently small 
approximations errors of $\mathbf{K}$ and $\mathbf{L}$ by \cite[Prop. A.2]{HumBek25b}, as well 
as sufficiently small approximation errors of $\mathbf{R}$ by the Lyapunov analysis in the proof of 
\cite[Prop. 2.1]{HuLVaz19} (see Footnote~\ref{fn:rstab} for the analogous argument in the 
$\infty+\infty$ case). However, the control law \eqref{eq:infUappr2} additionally contains 
persistent error terms due to the macro measurements.
In order to derive \eqref{eq:macrostab}, we view $\left(\begin{smallmatrix}
\tilde{u}(t) \\ \tilde{v}(t)\end{smallmatrix}\right) - \mathcal{F} \left(\begin{smallmatrix}
\mathbf{u}(t) \\ \mathbf{v}(t) \end{smallmatrix}\right)$ and $\tilde{u}(t,1,\cdot) - \mathcal{F}_n 
\mathbf{u}(t,1)$ as bounded perturbations due to \eqref{eq:uveps} and utilize input-to-state 
stability results for these perturbations. Hence, due to \eqref{eq:uveps}, the solution to the 
$n+m$ system under the proposed control law satisfies \eqref{eq:macrostab}, e.g., by \cite[Thm 
3.18, Def. 3.17, Rem. 3.14]{MirPri20}.
\end{proof}
\end{theorem}

In Proposition~\ref{thm:ave} below, we provide a specific case of Theorem~\ref{thm:infUappr2}, in 
which we further assume availability of only an average (over the ensemble variables) of the 
macro measurements, thus relaxing the requirement of availability of macro measurements for 
each value of the ensemble variables. Since on the way of proving Proposition~\ref{thm:ave} we 
establish that the continuum system can be stabilized using averaged, continuum 
measurements/kernels, such a setup may also be useful in the case in which the control objective 
is stabilization of the (macro) continuum system itself, rather than stabilization of the (micro) 
large-scale system.

However, in order for the average measurements to be accurate enough, so that the 
corresponding, closed-loop continuum system is exponentially stable, which is a prerequisite for 
the corresponding large-scale, $n+m$ system to be stable (or it is viewed as a standalone result 
when the purpose practically is stabilization of the continuum system), we need to assume that 
the parameters of the continuum system \eqref{eq:inf}, \eqref{eq:infbc} feature small variations 
with respect to the ensemble variables. This in turn translates to an assumption that all 
parameters of the corresponding large-scale, $n+m$ system \eqref{eq:nmm}, \eqref{eq:nmmbc} 
are close to each other (in a specific sense). Furthermore, since \eqref{eq:uveps} has to hold for 
$t \geq 0$, we also need to impose a (technical) uniform boundedness assumption on the 
autonomous system \eqref{eq:inf}, \eqref{eq:infbc}. We are now ready to state and prove this 
result.

\begin{proposition}
\label{thm:ave}
Consider an $n+m$ system \eqref{eq:nmm}, \eqref{eq:nmmbc} with 
parameters satisfying Assumption~\ref{ass:nm}, such that the solution to the autonomous system 
(i.e., \eqref{eq:nmm}, \eqref{eq:nmmbc} with $\mathbf{U}=0$) for any initial condition 
$\left(\begin{smallmatrix} 
\mathbf{u}_0 \\ \mathbf{v}_0 \end{smallmatrix}\right) \in E$ is uniformly bounded in time. 
Assume that there exist some 
$\bar{\lambda},\bar{\mu}\in C^1([0,1]; \mathbb{R})$ and $\bar{\sigma},\bar{\theta},\bar{w}\in 
C([0,1]; \mathbb{R})$\footnote{We 
tacitly take $\bar{\psi}=0$ due to Assumption~\ref{ass:nm}.}  and 
$\bar{r}, \bar{q} \in \mathbb{R}$ such that, for all $i,\ell \in \{1,\ldots,n\}$ and $j,p \in 
\{1,\dots,m\}$,\footnote{This is equivalent to assuming that the components of the parameters are 
close to each other in the sense of \eqref{eq:nmavep}.}
\begin{subequations}
\label{eq:nmavep}
\begin{align}
\max_{x\in[0,1]}| \lambda_i(x) - \bar{\lambda}(x) | + \max_{x\in[0,1]}| \lambda'_i(x) - 
\bar{\lambda}'(x) |& \leq \bar{\varepsilon}, \\
\max_{x\in[0,1]}| \mu_j(x) - \bar{\mu}(x) | + \max_{x\in[0,1]}| \mu'_j(x) - \bar{\mu}'(x) |& \leq 
\bar{\varepsilon}, \\
\max_{x\in[0,1]}| \sigma_{i,\ell}(x) - \bar{\sigma}(x) | & \leq 
\bar{\varepsilon}, \\
\max_{x\in[0,1]}| w_{i,p}(x) - \bar{w}(x) | & \leq 
\bar{\varepsilon}, \\
\max_{x\in[0,1]}| \theta_{j,\ell}(x) - \bar{\theta}(x) | & \leq 
\bar{\varepsilon}, \\
\max_{x\in[0,1]}| \psi_{j,p}(x)| & \leq 
\bar{\varepsilon},\\
|q_{i,p} - \bar{q}| \leq 
\bar{\varepsilon}, \qquad
| r_{j,\ell} - \bar{r} | & \leq 
\bar{\varepsilon},
\end{align}
\end{subequations}
for some $\bar{\varepsilon} > 0$ sufficiently small. Construct a respective 
continuum system \eqref{eq:inf}, \eqref{eq:infbc} under the 
conditions of Theorem~\ref{thm:infsolappr} with $U = \mathcal{F}_m\mathbf{U}$ and 
$\varepsilon$ in \eqref{eq:infappr} sufficiently small\footnote{For example, under 
\eqref{eq:infcap} and sufficiently large $n$ and $m$.}, such that the autonomous continuum 
system has a uniformly bounded solution.  Assume that the macro measurements are given by 
\begin{equation}
\label{eq:uvave}
\tilde{u}(t,x,y) \equiv \int\limits_0^1u(t,x,\zeta)d\zeta, \quad 
\tilde{v}(t,x,\eta) \equiv \int\limits_0^1v(t,x,\zeta)d\zeta,
\end{equation}
for all $t \geq 0$ and almost all $x \in [0,1]$. Then, the 
 $n+m$ system \eqref{eq:nmm}, \eqref{eq:nmmbc} under the control law 
\eqref{eq:infUappr2}--\eqref{eq:infUappr2k}, \eqref{eq:uvave} satisfies, for all $t \geq 0$,
\begin{equation}
\label{eq:dinfconv}
\resizebox{.99\columnwidth}{!}{$\displaystyle 
\left\| \left( \begin{smallmatrix}
\mathbf{u}(t) \\ \mathbf{v}(t)
\end{smallmatrix} \right) \right\|_E \leq M_c e^{-\omega_ct}\left(\left\| \left( \begin{smallmatrix}
\mathbf{u}_0 \\ \mathbf{v}_0
\end{smallmatrix} \right) \right\|_E + \varepsilon_{u,v} \right)+ \delta_\infty^1\left\| \left( 
\begin{smallmatrix}
\mathbf{u}_0 \\ \mathbf{v}_0
\end{smallmatrix} \right) \right\|_E + \delta_\infty^2$},
\end{equation}
for some positive constants $M_c, \omega_c, \delta_{\infty}^{1,2}$, such that 
$\delta_\infty^{1,2} \to 0$ as $\varepsilon, \varepsilon_{u,v} \to 0$.\footnote{For example, as $n, 
m \to\infty$; see Remark~\ref{rem:nminf}. Conceptually, estimate \eqref{eq:dinfconv} is expected, 
as, when viewing the continuum approximation error of the solutions as measurement error (in 
closed loop), by linearity and the uniform boundedness assumption of the open-loop system, one 
obtains a type of output-to-state stability property with respect to that measurement error. This 
error depends (via $\delta^2_\infty$) on the error due to continuum approximation of initial 
conditions (quantified by $\varepsilon_{u,v}$) and the error due to continuum parameters 
approximation (quantified by $\varepsilon$). The latter gives rise to an approximation error in the 
solutions operators that naturally grows with the size of initial conditions.}

\begin{proof}
Firstly, we show that the closed-loop system of \eqref{eq:inf}, \eqref{eq:infbc} 
under controls \eqref{eq:infUappr2}, \eqref{eq:uvave} is exponentially stable. Inserting 
\eqref{eq:uvave} to \eqref{eq:infUappr2}, we get 
\begin{align}
\label{eq:infUappr22}
\mathbf{U}(t) & = -\int\limits_{0}^1\bar{\mathbf{R}}u(t,1,\zeta)d\zeta + 
\int\limits_0^1\int\limits_0^1 \bar{\mathbf{K}}(1,\xi)u(t,\xi,\zeta)d\xi d\zeta \nonumber \\
& \qquad + \int\limits_0^1 \int\limits_0^1 \bar{\mathbf{L}}(1,\xi)v(t,\xi,\zeta)d\xi d\zeta,
\end{align}
where we denote $\bar{\mathbf{R}} = \int\limits_0^1 \widetilde{\mathbf{R}}(\zeta)d\zeta$, 
$\bar{\mathbf{K}}(1,\xi) = \int\limits_0^1 \mathbf{K}(1,\xi,\zeta)d\zeta$, and $\bar{\mathbf{L}}(1,\xi) 
= 
\int\limits_0^1 \mathbf{L}(1,\xi,\zeta)d\zeta$. Comparing \eqref{eq:infUappr22} with the exact 
backstepping control law 
\eqref{eq:infU}, we have\footnote{We interchangeably view 
$\bar{\mathbf{R}},\bar{\mathbf{K}},\bar{\mathbf{L}}$ as constant functions in $\zeta$.}
\begin{align}
\label{eq:infUcomp}
\mathcal{F}_m\mathbf{U}(t) & = U(t,\cdot) + \int\limits_0^1 \left(R(\cdot,\zeta) - 
\mathcal{F}_m\bar{\mathbf{R}}\right)u(t,1,\zeta)d\zeta \nonumber \\
& \quad - \resizebox{.75\columnwidth}{!}{$\displaystyle \int\limits_0^1\int\limits_0^1 
\left(K(1,\xi,\cdot,\zeta) - \mathcal{F}_m\bar{\mathbf{K}}(1,\xi)\right)u(t,\xi,\zeta)d\xi d\zeta$}
\nonumber \\
& \quad - \resizebox{.75\columnwidth}{!}{$\displaystyle\int\limits_0^1 \int\limits_0^1 
\left(L(1,\xi,\cdot,\zeta)  - 
\mathcal{F}_m\bar{\mathbf{L}}(1,\xi)v(t,\xi,\zeta)\right) d\xi 
d\zeta$}.
\end{align}
Since the exact backstepping control $U$ exponentially stabilizes the continuum system 
\eqref{eq:inf}, \eqref{eq:infbc}, due to analogous arguments to the proof of 
Theorem~\ref{thm:infUappr2}\footnote{\label{fn:rstab}In particular, the additional, remaining term 
affecting the boundary of the target system \eqref{eq:infts}, \eqref{eq:inftsbc} as $\beta(t,\eta,1) 
= \int\limits_0^1 \left(R(\eta,\zeta) - 
\mathcal{F}_m\bar{\mathbf{R}}\right)u(t,1,\zeta)d\zeta$ can be dominated in the derivative of the 
Lyapunov functional \eqref{eq:lyapd} by the term $-e^{-\delta}\|\alpha(t,1,\cdot)\|_{L^2}^2$ 
(recall $\alpha  \equiv u$), provided that 
$\left\| \int\limits_0^1  \left(R(\cdot,\zeta) - 
\mathcal{F}_m\bar{\mathbf{R}}\right)d\zeta\right\|_{L^2}^2 \leq \frac{e^{-2\delta}}{D(1)}$.}, the 
closed-loop 
system 
\eqref{eq:inf}, \eqref{eq:infbc} under 
controls \eqref{eq:infUappr2}, \eqref{eq:uvave} is exponentially stable provided that 
\begin{subequations}
\begin{align}
\| R(\cdot,\cdot) - \mathcal{F}_m\bar{\mathbf{R}}\|_{L^2}, \label{eq:Rcomp} \\
\esssup_{\xi \in [0,1] }\| K(1,\xi,\cdot,\cdot) - \mathcal{F}_m\bar{\mathbf{K}}(1,\xi) \|_{L^2}, 
\label{eq:Kcomp} \\
\esssup_{\xi \in [0,1] }\|  L(1,\xi,\cdot,\cdot)  - \mathcal{F}_m\bar{\mathbf{L}}(1,\xi) \|_{L^2}, 
\label{eq:Lcomp}
\end{align}
\end{subequations}
are sufficiently small. For \eqref{eq:Rcomp}, we can estimate
\begin{align}
\label{eq:Rbarest1}
\| R(\cdot,\cdot) - \mathcal{F}_m\bar{\mathbf{R}}\|_{L^2} & \leq \nonumber \\
\| R - R^{m,n}\|_{L^2} + \| R^{m,n} - \bar{r}\|_{L^2} +  \| \bar{r} - 
\mathcal{F}_m\bar{\mathbf{R}}\|_{L^2},
\end{align}
where $\| \bar{r} -\mathcal{F}_m\bar{\mathbf{R}}\|_{L^2}$ can be estimated by, recalling 
\eqref{eq:Rtilde},
\begin{align}
\label{eq:Rbarest2}
\|\bar{r} - \mathcal{F}_m\bar{\mathbf{R}}\|^2_{L^2} & = \nonumber \\ 
\sum_{j=1}^m\int\limits_{(j-1)/m}^{j/m}\int\limits_0^1 \left(\bar{r} - 
m\int\limits_{(j-1)/m}^{j/m}\int\limits_0^1 R(\chi,\zeta)d\zeta d\chi\right)^2d\zeta d\chi & = 
\nonumber \\
\sum_{j=1}^m \frac{1}{m} \left(m\int\limits_{(j-1)/m}^{j/m}\int\limits_0^1(\bar{r} -  R(\chi,\zeta)
)d\zeta d\chi\right)^2 & \leq \nonumber \\
\sum_{j=1}^m \int\limits_{(j-1)/m}^{j/m}\int\limits_0^1(\bar{r} -  R(\chi,\zeta)
)^2d\zeta d\chi & = \nonumber \\
\|\bar{r} - R\|_{L^2}^2,
\end{align}
where we used the Cauchy-Schwarz inequality. Thus, by \eqref{eq:infappr}, \eqref{eq:nmavep}, 
\eqref{eq:Rbarest1}, and \eqref{eq:Rbarest2}, we have $\| R(\cdot,\cdot) - 
\mathcal{F}_m\bar{\mathbf{R}}\|_{L^2} \leq 
2(\varepsilon + \bar{\varepsilon})$.

For estimating \eqref{eq:Kcomp}, \eqref{eq:Lcomp}, we first note that, by \eqref{eq:infappr}, 
\eqref{eq:nmavep}, the continuum parameters $\lambda,\mu,\sigma,W,\theta,\psi,Q$ are close to 
$\bar{\lambda},\bar{\mu},\bar{\sigma},\bar{\theta},\bar{w},\bar{\psi}=0, \bar{q}$, respectively, for 
all 
$x \in [0,1]$ and in the $L^2$ sense in $(\eta,\zeta)$.\footnote{We interchangeably view 
$\bar{\lambda},\bar{\mu},\bar{\sigma},\bar{\theta},\bar{w}, \bar{q}$ as constant functions  in 
$(\eta,\zeta)$.} We then introduce kernels $\bar{K}, \bar{L}$ 
that are the solution of the $2\times 2$ kernel equations  \cite[(18),
(19)]{DiMVaz13} for the parameters 
$\bar{\lambda},\bar{\mu},\bar{\sigma},\bar{\theta},\bar{w},\bar{\psi}=0,
\bar{q}$, i.e.,
\begin{subequations}
\label{eq:2x2k}
\begin{align}
\bar{\mu}(x) \bar{K}_x(x,\xi) - \bar{\lambda}(\xi)\bar{K}_{\xi}(x,\xi)
  - \bar{\lambda}_{\xi}(\xi)\bar{K}(x,\xi) & = \nonumber \\
  \bar{\sigma}(\xi)\bar{K}(x,\xi) + \bar{\theta}(\xi)\bar{L}(x,\xi), \\
\bar{\mu}(x)\bar{L}_x(x,\xi) + \bar{\mu}(\xi)\bar{L}_{\xi}(x,\xi) +
  \bar{\mu}_{\xi}(\xi)L(x,\xi) & = \nonumber \\
  \bar{w}(\xi)\bar{L}(x,\xi),
\end{align}
\end{subequations}
with boundary conditions
\begin{subequations}
\label{eq:2x2kbc}
\begin{align}
\bar{K}(x,x) & = - \frac{\bar{\theta}(x)}{\bar{\lambda}(x) +
               \bar{\mu}(x)}, \\
\bar{L}(x,0) & = \frac{1}{\bar{\mu}(0)}\bar{q}\bar{\lambda}(0)\bar{K}(x,0).
\end{align}
\end{subequations}
The solution 
$\bar{K}, \bar{L}$ to the $2\times 2$ kernel equations, which is well-posed by 
\cite[Sect. V]{DiMVaz13}, satisfies \eqref{eq:infk}, 
\eqref{eq:infkbc} for $\bar{\lambda},\bar{\mu},\bar{\sigma},\bar{\theta},\bar{w},\bar{\psi}=0, 
\bar{q}$, 
when interpreted as constant functions in $\eta,\zeta \in [0,1]$. Note that \eqref{eq:infkbc1}, 
\eqref{eq:infkbc4} become redundant when $\mu$ is $\eta$-invariant (and 
$\psi=0$)\footnote{The proper form of \eqref{eq:infkbc1} 
is \eqref{eq:infkbcseg1}, which is trivially satisfied for $\eta$-invariant $\mu$ when $\psi=0$. The 
artificial boundary condition \eqref{eq:infkbc4}, \eqref{eq:infl} can be assigned in the same form, 
so that it is trivially satisfied as well.}, 
and 
\eqref{eq:infkbc2}, \eqref{eq:infkbc3} hold for all $\eta, \zeta \in [0,1]$ (due to invariance in 
$\eta,\zeta$). Hence, 
due to well-posedness of the kernel equations \eqref{eq:2x2k}, \eqref{eq:2x2kbc} and 
\eqref{eq:infk}, \eqref{eq:infkbc} by \cite[Sect. V]{DiMVaz13} and Theorem~\ref{thm:infkwp}, 
respectively, and the proximity of the parameters 
$\lambda,\mu,\sigma,W,\theta,\psi,Q$ to 
$\bar{\lambda},\bar{\mu},\bar{\sigma},\bar{\theta},\bar{w},\bar{\psi}=0, \bar{q}$, respectively, 
there exists some 
$\delta_{\bar{\varepsilon}} > 0$ depending continuously on $\bar{\varepsilon}$ with 
$\delta_{\bar{\varepsilon}} \to 0$ as $\bar{\varepsilon}\to 0$ such that 
\begin{subequations}
\begin{align}
\esssup_{(x,\xi) \in \mathcal{T}} \| \bar{K}(x,\xi,\cdot,\cdot) - 
K(x,\xi,\cdot,\cdot)\|_{L^2([0,1]^2; \mathbb{R})} & \leq \delta_{\bar{\varepsilon}}, \\
\esssup_{(x,\xi) \in \mathcal{T}} \| \bar{L}(x,\xi,\cdot,\cdot) - 
L(x,\xi,\cdot,\cdot)\|_{L^2([0,1]^2; \mathbb{R})} & \leq \delta_{\bar{\varepsilon}}.
\end{align}
\end{subequations}
Now, estimating \eqref{eq:Kcomp}, \eqref{eq:Lcomp} similarly to \eqref{eq:Rbarest1}, 
\eqref{eq:Rbarest2} for almost every $\xi \in [0,1]$, we obtain that 
\begin{subequations}
\begin{align}
\esssup_{\xi \in [0,1] }\| K(1,\xi,\cdot,\cdot) - \mathcal{F}_m\bar{\mathbf{K}}(1,\xi) \|_{L^2} & \leq 
2(\delta_\varepsilon + \delta_{\bar{\varepsilon}}), \\
\esssup_{\xi \in [0,1] }\|  L(1,\xi,\cdot,\cdot)  - \mathcal{F}_m\bar{\mathbf{L}}(1,\xi) \|_{L^2} & \leq 
2(\delta_\varepsilon + \delta_{\bar{\varepsilon}}),
\end{align}
\end{subequations}
where $\delta_\varepsilon, \delta_{\bar{\varepsilon}} \to 0$ as $\varepsilon, \bar{\varepsilon} \to 
0$, which concludes that the closed-loop system of \eqref{eq:inf}, \eqref{eq:infbc} under control 
$U=\mathcal{F}_m \mathbf{U}$ with 
\eqref{eq:infUappr2}, \eqref{eq:uvave} is exponentially stable, provided that $\varepsilon, 
\bar{\varepsilon}$ are sufficiently small. In particular, the following holds for some $M_c, 
\omega_c > 0$
\begin{equation}
\label{eq:contest}
\left\| \left(\begin{smallmatrix}
u(t) \\ v(t)
\end{smallmatrix} \right) \right\|_{E_c^2} \leq M_c e^{-\omega_c t}\left\| 
\left(\begin{smallmatrix}
u_0 \\ v_0
\end{smallmatrix} \right) \right\|_{E_c^2},\quad t\geq0.
\end{equation}

Secondly, the system comprising \eqref{eq:nmm}, \eqref{eq:nmmbc} and \eqref{eq:inf}, 
\eqref{eq:infbc} with controls \eqref{eq:infUappr2}, \eqref{eq:uvave}, and $U=\mathcal{F}_m 
\mathbf{U}$ with (57), (64), respectively, has a cascade structure, 
and employing the notation of Theorem~\ref{thm:infsolappr}, the solution can be written as 
\begin{equation}
\label{eq:cassol}
\begin{bmatrix}
\left( \begin{smallmatrix}
u^n(t) \\ v^m(t)
\end{smallmatrix}\right) \\ 
\left( \begin{smallmatrix}
u(t) \\ v(t)
\end{smallmatrix}\right)
\end{bmatrix} = 
\begin{bmatrix}
\mathbb{T}_t^{n,m} & \Phi_t^{n,m}\mathbb{K}^m_t \\ 0 & \mathbb{T}_t + \Phi_t\mathbb{K}^m_t
\end{bmatrix}\begin{bmatrix}
\left( \begin{smallmatrix}
u_0^n \\ v_0^m
\end{smallmatrix}\right) \\ 
\left( \begin{smallmatrix}
u_0 \\ v_0
\end{smallmatrix}\right)
\end{bmatrix},
\end{equation}
where we employed $\mathcal{F}$ to transform the $n+m$ system to $E_c^2$ and 
$\mathbb{K}_t^m$ is such that $U(t) = \mathcal{F}_m\mathbf{U}(t) = \mathbb{K}_t^m\left( 
\begin{smallmatrix}
u_0 \\ v_0
\end{smallmatrix}\right)$ in the closed-loop system \eqref{eq:inf}, \eqref{eq:infbc},
\eqref{eq:infUappr22}. As $\mathbb{T}_t^{n,m}$ is bounded by 
assumption, and as we showed that $\mathbb{T}_t + \Phi_t\mathbb{K}^m_t$ 
is exponentially stable, the solution given by \eqref{eq:cassol} is bounded, provided that 
$\|\Phi_t^{n,m}\mathbb{K}_t^m\|_{\mathcal{L}(E_c^2)}$ is bounded uniformly in time, which 
follows 
by \cite[Prop. 4.3.3, Prop. 4.3.6, Prop. 4.4.5]{TucWeiBook} as
\begin{align}
\label{eq:phikest}
\|\Phi_t^{n,m}\mathbb{K}_t^m\|_{\mathcal{L}(E_c^2)} & \leq M_{\Phi^{n,m}}'\|e^{\omega' 
t}\mathbb{K}^m_t\|_{\mathcal{L}(E_c,L^2([0,1]; L^2([0,1]; \mathbb{R})))} \nonumber \\
&  \leq M_{\Phi^{n,m}}'M_{\mathbb{K}^m}
\end{align}
for some constants $M_{\Phi^{n,m}}', M_{\mathbb{K}^m } > 0$ and any $0 < \omega' < 
\omega_c$.

Finally, we can reuse Theorem~\ref{thm:infsolappr} and 
construct a continuum approximation of the 
cascade system by replacing the $n+m$ system \eqref{eq:nmm}, \eqref{eq:nmmbc} by its 
continuum approximation \eqref{eq:inf}, \eqref{eq:infbc} constructed in the statement of the 
theorem, so that the solution to the continuum approximation of the cascade system is given by
\begin{equation}
\label{eq:cassolappr}
\begin{bmatrix}
\left( \begin{smallmatrix}
u(t) \\ v(t)
\end{smallmatrix}\right) \\ 
\left( \begin{smallmatrix}
u(t) \\ v(t)
\end{smallmatrix}\right)
\end{bmatrix} = 
\begin{bmatrix}
\mathbb{T}_t & \Phi_t\mathbb{K}^m_t \\ 0 & \mathbb{T}_t + \Phi_t\mathbb{K}^m_t
\end{bmatrix}\begin{bmatrix}
\left( \begin{smallmatrix}
u_0 \\ v_0
\end{smallmatrix}\right) \\ 
\left( \begin{smallmatrix}
u_0 \\ v_0
\end{smallmatrix}\right)
\end{bmatrix}.
\end{equation}
The difference of \eqref{eq:cassol} and \eqref{eq:cassolappr} can be estimated as (we omit the 
second components as those are identical)
\begin{align}
\label{eq:contapprU0}
\left\|\left( \begin{smallmatrix}
u^n(t) \\ v^m(t)
\end{smallmatrix}\right) - \left( \begin{smallmatrix}
u(t) \\ v(t)
\end{smallmatrix}\right) \right\|_{E_c^2} 
\leq \nonumber \\
\left\| (\mathbb{T}_t^{n,m} - \mathbb{T}_t
)\left(\begin{smallmatrix}
u_0^n \\ v_0^m
\end{smallmatrix}\right) \right\|_{E_c^2} 
+ \|\mathbb{T}_t\|_{\mathcal{L}(E_c^2)} \left\| \left( \begin{smallmatrix}
u_0^n \\ v_0^m
\end{smallmatrix}\right)  - \left( \begin{smallmatrix}
u_0 \\ v_0
\end{smallmatrix}\right) \right\|_{E_c^2} \nonumber \\
+\|(\Phi_t^{n,m} - \Phi_t)\mathbb{K}_t^m \left(\begin{smallmatrix}
u_0 \\ v_0
\end{smallmatrix}\right) \|_{E_c^2},
\end{align}
where all the terms are uniformly bounded in time due to boundedness of $\mathbb{T}_t^{n,m}, 
\mathbb{T}_t$ (by assumption) and \eqref{eq:phikest}\footnote{An analogous estimate to  
\eqref{eq:phikest} holds for 
$\Phi_t\mathbb{K}_t^m$, because $\mathbb{T}_t$ is uniformly bounded by construction, enabled 
by the uniform boundedness assumption on $T_t^{n,m}$.}, and 
they tend to 
zero as $\varepsilon, \varepsilon_{u,v}\to 
0$ analogously to the proof of Theorem~\ref{thm:infsolappr}, where we can take $T=\infty$ due 
to both solutions being uniformly bounded in time. Thus, there exist some 
$\delta_{\mathbb{T}}, M_{\mathbb{T}}, \delta_{\Phi}> 0$, where $\delta_{\mathbb{T}},\delta_{\Phi} 
\to 0$ as $\varepsilon \to 0$ (due to continuous dependence of the solution operators to the 
respective parameters), such that
\begin{align}
\label{eq:dinfty}
\sup_{t \geq 0}\left\|\mathcal{F} \left( \begin{smallmatrix}
\mathbf{u}(t) \\ \mathbf{v}(t)
\end{smallmatrix}\right) - \left( \begin{smallmatrix}
u(t) \\ v(t)
\end{smallmatrix}\right) \right\|_{E_c^2} & \leq \left(\delta_{\mathbb{T}} +  
\delta_{\Phi}\right)\left\| \left(\begin{smallmatrix}
u_0^n \\ v_0^m
\end{smallmatrix}\right) \right\|_{E_c^2} \nonumber \\
& \quad + \left(M_{\mathbb{T}} +  
\delta_{\Phi}\right)\varepsilon_{u,v}.
\end{align}
Hence, \eqref{eq:dinfconv} follows with $\delta_\infty^1 = \delta_{\mathbb{T}}+\delta_{\Phi}$ and 
$\delta_\infty^2 = \left(M_{\mathbb{T}} +  
\delta_{\Phi}\right)\varepsilon_{u,v}$ by employing the triangle inequality
\begin{equation}
\label{eq:}
\left\| \left( \begin{smallmatrix}
u^n(t) \\ v^m(t)
\end{smallmatrix}\right) \right\|_{E_c^2} \leq \left\| \left( \begin{smallmatrix}
u(t) \\ v(t)
\end{smallmatrix}\right) \right\|_{E_c^2} + \left\| \left( \begin{smallmatrix}
u^n(t) \\ v^m(t)
\end{smallmatrix}\right) - \left( \begin{smallmatrix}
u(t) \\ v(t)
\end{smallmatrix}\right) \right\|_{E_c^2},
\end{equation}
together with \eqref{eq:contest}, \eqref{eq:contapprU0}, and \eqref{eq:dinfty}.
\end{proof}
\end{proposition}

\section{Continuum Kernels Well-Posedness} \label{sec:infkwp}

\begin{theorem}
\label{thm:infkwp}
Under Assumption~\ref{ass:inf}, the kernel equations \eqref{eq:infk}--\eqref{eq:infl} have a 
well-posed solution $K, L \in L^\infty(\mathcal{T}; L^2([0,1]^2; \mathbb{R}))$.
\end{theorem}

The proof is presented at the end of this section by utilizing the following lemmas.

\begin{lemma}[Splitting the kernel equations to subdomains]
The kernel equations \eqref{eq:infk} can be equivalently written in $L^\infty(\mathcal{T}; 
L^2([0,1]^2; \mathbb{R}))^2$ as 
\begin{subequations}
	\label{eq:infkseg}
	\begin{align}
		\mu(x,\eta)K^i_x(x,\xi,\eta,\zeta) - K^i_\xi(x,\xi,\eta,\zeta)\lambda(\xi,\zeta) & \nonumber \\ 
		- K^i(x,\xi,\eta,\zeta)\lambda_\xi(\xi,\zeta) 
		& = \nonumber \\
		\resizebox{.93\columnwidth}{!}{$\displaystyle 
		\int\limits_0^1K^i(x,\xi,\eta,\chi)\sigma(\xi,\chi,\zeta)d\chi +
		\int\limits_0^1L^i(x,\xi,\eta,\chi)\theta(\xi,\chi,\zeta)d\chi$}, \label{eq:infkseg1}
		& \\
		\mu(x,\eta)L^i_x(x,\xi,\eta,\zeta) + L^i_\xi(x,\xi,\eta,\zeta)\mu(\xi,\zeta) & \nonumber \\
		+ L^i(x,\xi,\eta,\zeta)\mu_\xi(\xi,\zeta) & = \nonumber \\
		\resizebox{.93\columnwidth}{!}{$\displaystyle 
		\int\limits_0^1K^i(x,\xi,\eta,\chi)W(\xi,\chi,\zeta)d\chi + 
		\int\limits_0^1L^i(x,\xi,\eta,\chi)\psi(\xi,\chi,\zeta)d\chi$},  \label{eq:infkseg2}
		&
	\end{align}	
\end{subequations}
for $i\in \{a,b,c\}$, where $K^i, L^i$ denote the restrictions of the kernels to $\mathcal{H}_i$ 
defined as
\begin{subequations}
\label{eq:Hseg}	
\begin{align}
\mathcal{H}_a & = \left\{(x,\xi,\eta,\zeta) \in [0,1]^4: \eta \leq \zeta, \xi(\eta,\zeta) \leq 
\rho(x,\eta,\zeta)\right\}, \\
\mathcal{H}_b & = \resizebox{.87\columnwidth}{!}{$\displaystyle \left\{(x,\xi,\eta,\zeta) \in 
[0,1]^4: 
\eta 
\leq \zeta,
\rho(x,\eta,\zeta)\leq \xi(\eta,\zeta) \leq x \right\}$}, \\
\mathcal{H}_c & = \left\{(x,\xi,\eta,\zeta) \in [0,1]^4: \zeta < \eta, \xi \leq x \right\},
\end{align}	
\end{subequations}
where\footnote{Note that $\xi(\eta,\zeta)= \rho(x,\eta,\zeta)$ is the characteristic hypersurface 
of 
\eqref{eq:infk2}.}
\begin{equation}
\rho(x,\eta,\zeta) = \phi_\zeta^{-1}\left(\phi_\eta(x)\right),
\end{equation}
for $x\in[0,1]$ and $0 \leq \eta \leq \zeta \leq 1$ with $\phi_\eta$ (respectively $\phi_\zeta$) 
given 
by
\begin{equation}
\phi_\eta(x) = \int\limits_0^x \frac{ds}{\mu(s,\eta)}.
\end{equation}
The boundary conditions for \eqref{eq:infkseg} are given by
\begin{subequations}
	\label{eq:infkbcseg}
	\begin{align}
		\psi(x,\eta,\zeta) &  = \mu(x,\zeta)L^j(x,x,\eta,\zeta) - L^j(x,x,\eta,\zeta)\mu(x,\eta), 
		\label{eq:infkbcseg1} \\
		-\theta(x,\eta,\zeta) & =  \mu(x,\eta)K^j(x,x,\eta,\zeta) + K^j(x,x,\eta,\zeta)\lambda(x,\zeta), 
		\label{eq:infkbcseg2} \\
		 L^a(x,0,\eta,\zeta) & = 
		\frac{1}{\mu(0,\zeta)}\int\limits_0^1K^a(x,0,\eta,\chi)\lambda(0,\chi)Q(\chi,\zeta)d\chi, 
		\label{eq:infkbcseg3} \\
		L^c(1,\xi,\eta,\zeta) & = l(\xi,\eta,\zeta), \label{eq:infkbcseg4}
	\end{align}	
\end{subequations}
where $j \in \{b,c\}$, in addition to which the $K$ kernel is subject to the continuity condition
\begin{equation}
 \label{eq:infkbcsegc}
K^a(x,\rho(x,\eta,\zeta),\eta,\zeta) = K^b(x,\rho(x,\eta,\zeta),\eta,\zeta).
\end{equation}

\begin{proof}
The backstepping transformation \eqref{eq:infV2} can be written in terms of the segmented 
kernels as
\begin{align}
\label{eq:infVseg}
\beta(t,x,\eta) & = v(t,x,\eta) - \int\limits_\eta^1 \int\limits_0^{\rho(x,\eta,\zeta)}
L^a(x,\xi,\eta,\zeta)v(t,\xi,\zeta)d\xi d\zeta \nonumber \\
& \qquad  - \int\limits_\eta^1 \int\limits_{\rho(x,\eta,\zeta)}^xL^b(x,\xi,\eta,\zeta)v(t,\xi,\zeta)d\xi 
d\zeta \nonumber 
\\
& \qquad - \int\limits_0^\eta \int\limits_0^x L^c(x,\xi,\eta,\zeta)v(t,\xi,\zeta)d\xi d\zeta 
\nonumber \\
& \qquad - \int\limits_\eta^1 \int\limits_0^{\rho(x,\eta,\zeta)}
K^a(x,\xi,\eta,\zeta)u(t,\xi,\zeta)d\xi d\zeta \nonumber \\
& \qquad  - \int\limits_\eta^1 \int\limits_{\rho(x,\eta,\zeta)}^xK^b(x,\xi,\eta,\zeta)u(t,\xi,\zeta)d\xi 
d\zeta \nonumber 
\\
& \qquad - \int\limits_0^\eta \int\limits_0^x K^c(x,\xi,\eta,\zeta)u(t,\xi,\zeta)d\xi d\zeta.
\end{align}
The segmented kernel equations \eqref{eq:infkseg} are obtained by inserting \eqref{eq:infVseg} to 
\eqref{eq:infts2} and integrating by parts once as in Appendix~\ref{app:ker}. In fact, the kernel 
equations \eqref{eq:infkseg} and boundary conditions \eqref{eq:infkbcseg} are of the same form 
as the ones presented in \eqref{eq:infk}, \eqref{eq:infkbc}, with the addition of the continuity 
condition \eqref{eq:infkbcsegc} that arises due to the segmentation of the domain $\mathcal{T} 
\times [0,1]^2$ when differentiating \eqref{eq:infVseg} in $x$ and integrating by parts once.
\end{proof}
\end{lemma}

\begin{figure}[!htb]
\begin{center}
\includegraphics[scale=.69]{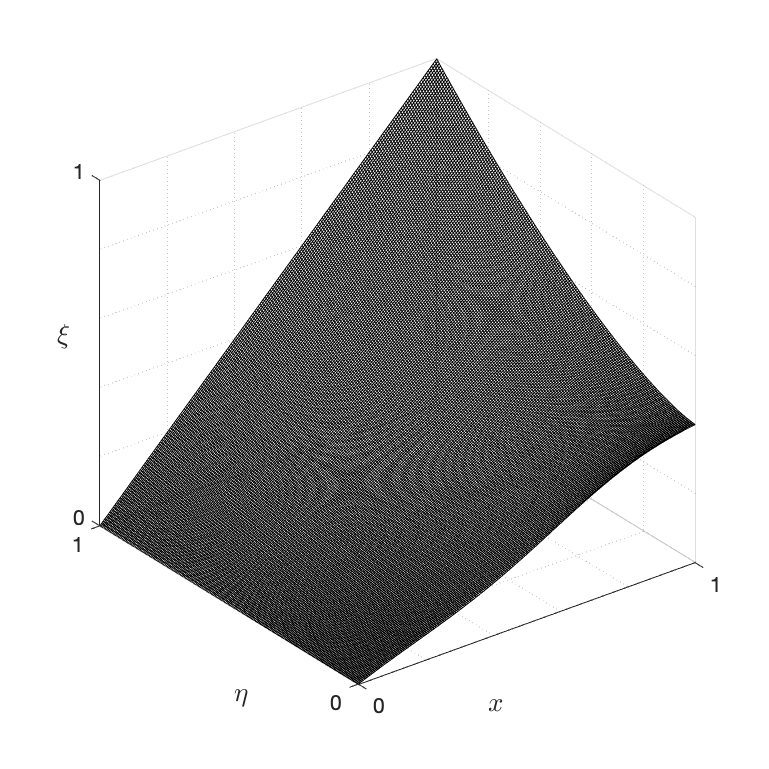}
\end{center}
\caption{Illustration of $\xi = \rho(x,\eta,\zeta)$ projected on $\zeta=1$. The characteristic 
hypersurface is a collection of such surfaces over all $0 \leq \eta \leq \zeta \leq 1$, which all 
contain the line $\xi = x$ at $\eta=\zeta$.}
\label{fig:seg}
\end{figure}

\begin{remark}
Note that the potential kernel discontinuities may 
only occur in $L$ for $\eta \leq \zeta$ on the hypersurface $\xi = \rho(x,\eta,\zeta)$, which is 
continuous 
and monotonic in all variables; see also Remark~\ref{rem:infk}.\footnote{Compared 
to the case of finite $m$, the 
characteristic hypersurface can be viewed as an infinite collection of characteristic curves for 
the 
$n+m$ kernels or $\infty+m$ kernels (see, e.g., \cite{HuLVaz19, HumBek26}). In fact, for 
any fixed $\zeta$ and $\eta$  such that $\eta \leq \zeta$, the characteristic 
hypersurface reduces to such a characteristic curve.} An illustration of the 
characteristic hypersurface projected on $\zeta=1$ is provided in Fig.~\ref{fig:seg}.
\end{remark}

\begin{lemma}[Continuity of characteristic projections]
The characteristic projections of the kernel equations \eqref{eq:infkseg} are continuous on 
$\mathcal{H}_i$ for $i \in \{a,b,c\}$.

\begin{proof}
Since $\mu, \lambda \in C^1([0,1]^2;\mathbb{R})$ by Assumption~\ref{ass:inf}, we can argue 
pointwise in $\eta,\zeta \in [0,1]$ and solve the characteristic projections for the $K^i$ and 
$L^i$ 
kernels for $i \in \{a,b,c\}$. For fixed (albeit arbitrary) $\eta,\zeta \in [0,1]$, the characteristic 
projections are, in fact, analogous to those encountered in the $n+m$ case, evolving in (a 
subset 
of) $(x,\xi) \in \mathcal{T}$. The characteristic projections for the $K^i$ kernels satisfy 
the following Cauchy problem on $s \in [0,s_f^i(\eta,\zeta)]$ for arbitrary, fixed $\eta,\zeta \in 
[0,1]$
\begin{subequations}
	\label{eq:infcpk}%
	\begin{align}
		\frac{d}{ds}\hat{x}^i(s,\eta,\zeta) & = -\mu\left(\hat{x}^i(s,\eta,\zeta),\eta\right), \\
		\frac{d}{ds}\hat{\xi}^i(s,\eta,\zeta) & = \lambda\left(\hat{\xi}^i(s,\eta,\zeta),\zeta\right),
	\end{align}
\end{subequations}
with boundary conditions $\hat{x}^i(0,\eta,\zeta) = x$, $\hat{x}^i(s_f^i,\eta,\zeta) = 
\hat{x}^i_f(\eta,\zeta)$,$\hat{\xi}^i(0,\eta,\zeta) = \xi$, $\hat{\xi}^i(s_f^i,\eta,\zeta) = 
\hat{\xi}^i_f(\eta,\zeta)$. 
Since $\mu(\cdot,\eta), \lambda(\cdot,\zeta)$ are continuously differentiable and positive by 
Assumption~\ref{ass:inf}, \eqref{eq:infcpk} has a unique continuously differentiable solution for 
any $(x,\xi) \in \mathcal{T}$ 
and for each $\eta,\zeta \in [0,1]$ by Picard---Lindel\"of theorem \cite[Thm 2.2]{TesBook}, 
where 
$\hat{x}^i$ is strictly decreasing in $s$ and $\hat{\xi}^i$ is strictly increasing in $s$. Thus, for 
$i = 
a$, the solution to \eqref{eq:infcpk} terminates at $s_f^a(\eta,\zeta)$ on 
$\hat{\xi}_f^a(\eta,\zeta) 
= \rho\left(\hat{x}_f^a(\eta,\zeta),\eta,\zeta\right)$, and the corresponding boundary condition 
is 
given by \eqref{eq:infkbcsegc}. For $i \in \{b,c\}$, the solution to \eqref{eq:infcpk} terminates at 
$s_f^i(\eta,\zeta)$ on $\hat{\xi}_f^i(\eta,\zeta) = \hat{x}_f^i(\eta,\zeta)$, and the corresponding 
boundary condition is given by \eqref{eq:infkbcseg2}.

Analogously to the characteristic projections for the $K^i$ kernels, we argue pointwise in 
$\eta,\zeta \in [0,1]$ to establish the characteristic projections for the $L^i$ kernels. For 
arbitrary, 
fixed $\eta,\zeta \in [0,1]$, the characteristic projections for the $L^i$ kernels satisfy the 
following Cauchy problem on $s \in [0, s_F^i(\eta,\zeta)]$
\begin{subequations}
	\label{eq:infcpl}%
	\begin{align}
		\frac{d}{ds}\hat{\chi}(s,\eta,\zeta) & = \epsilon(\eta,\zeta)\mu(\hat{\chi}(s,\eta,\zeta),\eta), \\
		\frac{d}{ds}\hat{\zeta}(s,\eta,\zeta) & = 
		\epsilon(\eta,\zeta)\mu(\hat{\zeta}(s,\eta,\zeta),\zeta),
	\end{align}
\end{subequations}
with boundary conditions $\hat{\chi}^i(0,\eta,\zeta) = x$, $\hat{\chi}^i(s_F^i,\eta,\zeta) = 
\hat{\chi}^i_F(\eta,\zeta)$, $\hat{\zeta}^i(0,\eta,\zeta) = \xi$, $\hat{\zeta}^i(s_F^i,\eta,\zeta) = 
\hat{\zeta}^i_F(\eta,\zeta)$, and with
\begin{equation}
\label{eq:eps}
	\epsilon(\eta,\zeta) = \begin{cases}
		1, & \eta > \zeta \\
		-1, & \eta \leq \zeta
	\end{cases}.
\end{equation}
The location of the terminal condition $\left(\hat{\chi}^i_F(\eta,\zeta), 
\hat{\zeta}_F^i(\eta,\zeta)\right)$ 
depends on $i \in \{a,b,c\}$ as follows.
\begin{itemize}
\item For $i = a$, we have $\eta \leq \zeta$, and hence, both $\hat{\chi}^a$ and $\hat{\zeta}^a$ 
are strictly decreasing in $s$.  Thus, the solution to \eqref{eq:infcpl} terminates at 
$s_F^a(\eta,\zeta)$ on 
$\hat{\zeta}_F^a(\eta,\zeta) = 0$, and the corresponding boundary condition is given by 
\eqref{eq:infkbcseg3}.
\item For $i = b$, both $\hat{\chi}^b$ and $\hat{\zeta}^b$ are also strictly decreasing in $s$, so 
that the solution to \eqref{eq:infcpl} terminates at $s_F^b(\eta,\zeta)$ on 
$\hat{\zeta}_F^b(\eta,\zeta) = \hat{\chi}_F^b(\eta,\zeta)$, and the corresponding boundary 
condition is given by \eqref{eq:infkbcseg1}.\footnote{Note that this only applies for $\eta < 
\zeta$, 
whereas for $\eta = \zeta$, the solution to \eqref{eq:infcpl} is parallel to $\xi=x$ and terminates 
on 
$\xi=0$, which is covered by the case $i=a$.} 
\item For $i=c$, both $\hat{\chi}^c$ and $\hat{\zeta}^c$ are strictly increasing in $s$, and the 
solution to \eqref{eq:infcpl} terminates at $s_F^c(\eta,\zeta)$ either on 
$\hat{\zeta}_F^c(\eta,\zeta) = \hat{\chi}_F^c(\eta,\zeta)$ or on $\hat{\chi}_F^c(\eta,\zeta) = 1$. 
The corresponding boundary condition is given either by \eqref{eq:infkbcseg1} or by 
\eqref{eq:infkbcseg4}, respectively.
\end{itemize}

In order to argue continuity of characteristic projections, we first need the 
mappings\footnote{Here we need to account for the dependence of $s_f^i, s_F^i$ on $(x,\xi) 
\in 
\mathcal{T}$, whereas in the above, $(x,\xi) \in \mathcal{T}$ was considered fixed, and hence, 
it 
was omitted (cf. \cite[Sect. VI.B]{AllKrs25}).} 
$(\eta,\zeta,x,\xi) \in [0,1]^2 \times \mathcal{T} \to s_f^i(\eta,\zeta,x,\xi)$ and $(\eta,\zeta,x,\xi) 
\in 
[0,1]^2 \times \mathcal{T} \to s_F^i(\eta,\zeta,x,\xi)$  to be 
Lipschitz continuous in $\mathcal{H}_i$ for all $i \in \{a,b,c\}$. This follows by \cite[Lem. 
4]{AllKrs25}, as we can analogously prove that the 
above mappings are Lipschitz independently in $\eta$ and $\zeta$ (for arbitrary, fixed $\zeta$ 
and $\eta$, respectively), which then implies Lipschitzness in $(\eta,\zeta,x,\xi)$ for the full 
mapping. Consequently, the characteristic curves are continuous by \cite[Cor. 1]{AllKrs25}.
\end{proof}
\end{lemma}

Integrating \eqref{eq:infkseg} along the characteristic projections and plugging in the boundary 
conditions \eqref{eq:infkbcseg}, \eqref{eq:infkbcsegc} gives (pointwise in $(x,\xi,\cdot,\cdot) \in 
\mathcal{H}_i$ and in the $L^2$ sense in $(\eta,\zeta)$)\footnote{We drop $(\eta,\zeta)$ from 
$\hat{x}^i(s,\eta,\zeta)$ and $\hat{\xi}^i(s,\eta,\zeta)$ for notational brevity.}
\begin{subequations}
\label{eq:infkie}%
\begin{align}
	K^i\left(x,\xi,\eta,\zeta\right) -
	B^i_1\left(x^i_f\left(\eta,\zeta\right),\eta,\zeta\right) & = \nonumber \\
	-\int\limits_0^{s_f^i(\eta,\zeta)}\left(
	K^i\left(\hat{x}^i\left(s\right),\hat{\xi}^i\left(s\right),\eta,\zeta\right)
	\lambda_{\xi}\left(\hat{\xi}^i\left(s\right),\zeta\right) \right.
	& \nonumber \\
	+
	\int\limits_0^1\left(K^i\left(\hat{x}^i\left(s\right),\hat{\xi}^i\left(s\right),\eta,\chi\right)
	\sigma\left(\hat{\xi}^i\left(s\right),\chi,\zeta\right)\right. &
	\nonumber \\ \left.\left.
	+
	L^i\left(\hat{x}^i\left(s\right),\hat{\xi}^i\left(s\right),\eta,\chi\right)
	\theta\left(\hat{\xi}^i\left(s\right),\chi,\zeta\right)\right)d\chi
	\right)ds,
	& \label{eq:infkie1} \\
	L^i\left(x,\xi,\eta,\zeta\right) -
	B_2^i\left(\hat{\star}^i\left(s_F^i(\eta,\zeta),\eta,\zeta\right)\right)
	& = \nonumber \\ 
	\epsilon(\eta,\zeta)\int\limits_0^{s_F^i(\eta,\zeta)} \left(
	L^i\left(\hat{x}^i(s),\hat{\xi}^i\left(s\right),\eta,\zeta\right)\mu_\xi\left(\hat{\xi}^i\left(s\right),\zeta\right)
	\right.  & \nonumber \\
	-  \int\limits_0^1 
	\left(K^i\left(\hat{x}^i\left(s\right),\hat{\xi}^i\left(s\right),\eta,\chi\right)
	W\left(\hat{\xi}^i\left(s\right),\chi,\zeta\right) \right. & \nonumber \\
	\left. \left.
	+L^i\left(\hat{x}^i\left(s\right),\hat{\xi}^i(s),\eta,\chi\right)
	\psi\left(\hat{\xi}^i(s),\chi,\zeta\right)\right)d\chi
	\right)ds, \label{eq:infkie2}
	\end{align}
\end{subequations} where, for $i \in \{a,b,c\}$,
\begin{subequations}
	\label{eq:infbctmp}%
	\begin{align} B_1^i(x,\eta,\zeta) & = \begin{cases}
	K^b(x,\rho(x,\eta,\zeta),\eta,\zeta), & i = a \\
 	-\frac{\theta(x,\eta,\zeta)}{\lambda(x,\eta) + \mu(x,\zeta)}, & i \in \{b,c\}
		\end{cases}, \label{eq:bctmp1} \\
	B_2^i(\star,\eta,\zeta) & =
	\begin{cases}
	\frac{1}{\mu(0,\zeta)}\resizebox{.38\columnwidth}{!}{$\displaystyle \int\limits_0^1 
	K^a(x,0,\eta,\chi)\lambda(0,\chi)Q(\chi,\zeta)d\chi$}, & i = a	\\
	\frac{\psi(x,\eta,\zeta)}{\mu(x,\zeta) - \mu(x,\eta)}, & i \in \{b,c\} \\
	l(\xi,\eta,\zeta), & i = c
	\end{cases},
\end{align}
\end{subequations}
denote the boundary conditions according to the terminal conditions of the characteristic 
projections.\footnote{In $B_2^i$, $\star$ refers to $x$ or $\xi$ depending on which 
boundary condition is applied.} In the next lemma, we establish well-posedness of the integral 
form \eqref{eq:infkie} of 
the kernel equations using successive approximations.

\begin{lemma}[Convergence of successive approximations]
\label{lem:sa}
For $i \in \{a,b,c\}$, denote by $\left(K_\ell^i\right)_{\ell=0}^\infty$ and 
$\left(L_\ell^i\right)_{\ell=0}^\infty$ the sequences of successive approximations for respective 
kernels $K^i, L^i$ in \eqref{eq:infkie}, \eqref{eq:infbctmp}, where we initialize $K_0^i, L_0^i$ 
to 
zero. Then, the sequences of successive approximations converge such that\footnote{We 
tacitly 
extend the segmented kernels by zero functions outside their 
respective domain $\mathcal{H}_i$, so that the $L^2$ norm over $(\eta,\zeta) \in [0,1]^2$ is 
well-defined.}
\begin{subequations}
\label{eq:infksac}
\begin{align}
\lim_{\ell\to\infty} \| K_\ell^i(x,\xi,\cdot,\cdot) - K^i(x,\xi,\cdot,\cdot)\|_{L^2} & = 0, \\
\lim_{\ell\to\infty} \| L_\ell^i(x,\xi,\cdot,\cdot) - L^i(x,\xi,\cdot,\cdot)\|_{L^2} & = 0,
\end{align}
\end{subequations}
for all $(x,\xi, \cdot, \cdot) \in \mathcal{H}_i$.

\begin{proof}
Denoting $\Delta K_\ell^i = K_{\ell+1}^i - K_\ell^i$ and $\Delta L_\ell^i 
= L_{\ell+1}^i - L_{\ell}^i$, we can write
\begin{equation}
\label{eq:diffKLinf}
K_\ell^i = \sum_{l=0}^\ell \Delta K_l^i, \qquad L_\ell^i = \sum_{l=0}^\ell \Delta L_l^i,
\end{equation}
due to the initialization $K_0^i = L_0^i = 0$. Hence, the convergence of the sequences of 
successive approximations is equivalent to the convergence of the series \eqref{eq:diffKLinf}, 
which follows by showing the following relations
\begin{subequations}
\label{eq:infsaest}
\begin{align}
\|\Delta K_\ell^i(x,\xi,\cdot,\cdot)\|_{L^2} & \leq 
M\frac{(M_{K,L}m_\Phi^{-1}M_\Phi)^\ell}{\ell!}, \\
\|\Delta L_\ell^i(x,\xi,\cdot,\cdot)\|_{L^2} & \leq 
M\frac{(M_{K,L}m_\Phi^{-1}M_\Phi)^\ell}{\ell!},
\end{align}	
\end{subequations}
for all $(x,\xi,\cdot,\cdot) \in \mathcal{H}_i$, where the coefficients are given by
\begin{subequations}
\begin{align}
M & = M_B + (1+M_Q^1) \max_{x\in 
[0,1]}\frac{\|\theta(x,\cdot,\cdot)\|_{L^2}}{m_\lambda + 
m_\mu}, \label{eq:seM} \\
M_{K,L} & = 2(1+M_Q^1)(M_\lambda^1 + M_\sigma + M_\theta) \nonumber \\
& \qquad  + 2(M_\mu^1 + M_W + M_\psi), 
\label{eq:seMKL}
\end{align}
\end{subequations}
where $m_\lambda, m_\mu$ are given in \eqref{eq:mlam}, \eqref{eq:mmu}, and
\begin{subequations}
\label{eq:MB}	
\begin{align}
M_B & = \max_{x\in[0,1]} \left(\int\limits_0^1\int\limits_0^1 
\left(\frac{\psi(x,\eta,\zeta)}{\mu(x,\eta) 
- \mu(x,\zeta)}\right)^2d\eta d\zeta\right)^\frac{1}{2}, \\
M_\Phi & = \resizebox{.83\columnwidth}{!}{$\displaystyle \max_{(x,\xi) \in \mathcal{T}} 
\left(\int\limits_0^1\int\limits_0^1 \left(
e^{xe^{-\gamma\epsilon(\eta,\zeta)} } - e^{\xi e^{\gamma\epsilon(\eta,\zeta)} } + 
e^{e^{\gamma}} 
\right)^2d\eta 
d\zeta\right)^\frac{1}{2}$},
\end{align}
\end{subequations}
where $\epsilon$ is given in \eqref{eq:eps} and $\gamma > 0$ is sufficiently large such that 
$\displaystyle \frac{M_\mu}{m_\mu} < e^{2\gamma - e^{-\gamma}}$, 
where $M_\mu = \displaystyle \max_{x,\eta \in [0,1]}\mu(x,\eta)$, and $m_\Phi > 0$ is 
sufficiently 
small 
such that
\begin{equation}
\label{eq:mphi}
m_{\Phi} < \min\left\{m_\mu e^{\gamma} - M_\mu e^{e^{-\gamma}-\gamma}, (m_\mu + 
m_\lambda)e^{-\gamma} \right\},
\end{equation}
and
\begin{subequations}
\begin{align}
M_\lambda^1 & = \max_{x,y \in [0,1]} \lambda_x(x,y), \quad M_\mu^1 = \max_{x,\eta \in [0,1]} 
\mu_x(x,\eta), \\
M_\sigma & = \max_{x\in [0,1]}\left\|\int\limits_0^1\sigma(x,\eta,\cdot)d\eta\right\|_{L^2}, \\
M_\theta & = \max_{x\in [0,1]}\left\|\int\limits_0^1\theta(x,\eta,\cdot)d\eta\right\|_{L^2}, \\
M_W & = \max_{x\in [0,1]}\left\|\int\limits_0^1 W(x,\eta,\cdot)d\eta\right\|_{L^2}, \\
M_\psi & = \max_{x\in [0,1]}\left\|\int\limits_0^1\psi(x,\eta,\cdot)d\eta\right\|_{L^2}, \\
M_Q^1 & = \max_{\eta,\zeta \in [0,1]} \frac{\lambda(0,\eta)}{\mu(0,\zeta)}\left\|\int\limits_0^1 
Q(\chi,\cdot)d\chi\right\|_{L^2}.
\end{align}
\end{subequations}

The key in deriving the estimates \eqref{eq:infsaest} is to show by induction that the 
nonnegative 
function 
$\Phi \in C(\mathcal{T}; L^\infty([0,1]^2; \mathbb{R}))$ given by 
\begin{equation}
\Phi(x,\xi,\eta,\zeta) = e^{xe^{-\gamma\epsilon(\eta,\zeta)} } - e^{\xi 
e^{\gamma\epsilon(\eta,\zeta)} } + 
e^{e^{\gamma}},
\end{equation}
satisfies
\begin{subequations}
\label{eq:infsaestpw}
\begin{align}
|\Delta K_\ell^i(x,\xi,\eta,\zeta)| & \leq 
M\frac{(M_{K,L}m_\Phi^{-1})^\ell\Phi(x,\xi,\eta,\zeta)^\ell}{\ell!}, \\
|\Delta L_\ell^i(x,\xi,\eta,\zeta)| & \leq 
M\frac{(M_{K,L}m_\Phi^{-1})^\ell\Phi(x,\xi,\eta,\zeta)^\ell}{\ell!},
\end{align}
\end{subequations}
for all $(x,\xi,\cdot,\cdot) \in \mathcal{H}_i$ and almost every $\eta,\zeta \in [0,1]$. For the 
induction step, we show the following inequalities
\begin{subequations}
\label{eq:MPhi}
\begin{align}
\resizebox{.98\columnwidth}{!}{$\displaystyle\int\limits_0^{s_f^i(\eta,\zeta)} 
\Phi(\hat{x}^i(s,\eta,\zeta), 
\hat{\xi}^i(s,\eta,\zeta),\eta,\zeta)^\ell ds \leq 
\frac{1}{m_\Phi}
\frac{\Phi(x,\xi,\eta,\zeta)^{\ell+1}}{\ell+1}$}, \label{eq:MPhi1} \\
\resizebox{.98\columnwidth}{!}{$\displaystyle \int\limits_0^{s_F^i(\eta,\zeta)} 
\Phi(\hat{\chi}^i(s,\eta,\zeta), \hat{\zeta}^i(s,\eta,\zeta),\eta,\zeta)^\ell ds \leq 
\frac{1}{m_\Phi} \frac{\Phi(x,\xi,\eta,\zeta)^{\ell+1}}{\ell+1}$}, \label{eq:MPhi2}
\end{align}
\end{subequations}
for all $(x,\xi, \cdot,\cdot) \in \mathcal{H}_i$. We introduce 
a change of variables in \eqref{eq:MPhi1} as $\tau(s) = f_{i,\eta,\zeta}(s)$, where
\begin{subequations}
\begin{align}
\tau(s) & = \Phi(\hat{x}^i(s,\eta,\zeta), \hat{\xi}^i(s,\eta,\zeta),\eta,\zeta),
\end{align}
with
\begin{align}
d\tau & = 
-\left(e^{-\gamma\epsilon(\eta,\zeta)}e^{\hat{x}^i(s,\eta,\zeta)e^{-\gamma\epsilon(\eta,\zeta)}}\mu(\hat{x}^i
(s,\eta,\zeta),\eta) \right. \nonumber  \\
& \qquad \left. + 
e^{\gamma\epsilon(\eta,\zeta)}e^{\hat{\xi}^i(s,\eta,\zeta)e^{\gamma\epsilon(\eta,\zeta)}}\lambda(\hat{\xi}^i
(s,\eta,\zeta),\zeta)\right)ds \nonumber \\
& =: f(\hat{x}^i(s,\eta,\zeta), \hat{\xi}^i(s,\eta,\zeta),\eta,\zeta)ds,
\end{align}
\end{subequations}
so that \eqref{eq:MPhi1} becomes, denoting $\bar{x}^i(\tau,\eta,\zeta) = 
\hat{x}^i\left(f_{i,\eta,\zeta}^{-1}(s),\eta,\zeta\right)$ and $\bar{\xi}^i(\tau,\eta,\zeta) = 
\hat{\xi}^i\left(f_{i,\eta,\zeta}^{-1}(s),\eta,\zeta\right)$,
\begin{align}
\int\limits_0^{s_f^i(\eta,\zeta)} \Phi(\hat{x}^i(s,\eta,\zeta), \hat{\xi}^i(s,\eta,\zeta),\eta,\zeta)^\ell 
ds 
& =
\nonumber 
\\
\int\limits_{\Phi(x,\xi,\eta,\zeta)}^{\Phi\left(x_f^i(\eta,\zeta),\xi_f^i(\eta,\zeta),\eta,\zeta\right)} 
\frac{\tau^\ell 
d\tau}{f(\bar{x}^i(\tau,\eta,\zeta), \bar{\xi}^i(\tau,\eta,\zeta),\eta,\zeta)} & \leq \nonumber \\
\frac{1}{(m_\mu+m_\lambda)e^{-\gamma}}\frac{\Phi(x,\xi,\eta,\zeta)^{\ell+1}}{\ell+1}.
\end{align}
Similarly, we introduce a change of variables in \eqref{eq:MPhi2} as $\tau(s) = 
g_{i,\eta,\zeta}(s)$, 
where
\begin{subequations}
\begin{align}
\tau(s) & = \Phi(\hat{\chi}^i(s,\eta,\zeta), \hat{\zeta}^i(s,\eta,\zeta),\eta,\zeta),
\end{align}
with
\begin{align}
d\tau & = 
\epsilon(\eta,\zeta)\left(e^{-\gamma\epsilon(\eta,\zeta)}e^{\hat{\chi}^i(s,\eta,\zeta)e^{-\gamma\epsilon(\eta,\zeta)}}\mu(\hat{\chi}^i
(s,\eta,\zeta),\eta) \right. \nonumber  \\
& \qquad \left. - 
e^{\gamma\epsilon(\eta,\zeta)}e^{\hat{\zeta}^i(s,\eta,\zeta)e^{\gamma\epsilon(\eta,\zeta)}}\mu(\hat{\zeta}^i
(s,\eta,\zeta),\zeta)\right)ds \nonumber \\
& =: g(\hat{\chi}^i(s,\eta,\zeta), \hat{\zeta}^i(s,\eta,\zeta),\eta,\zeta)ds,
\end{align}
\end{subequations}
so that \eqref{eq:MPhi2} becomes, denoting $\bar{\chi}^i(\tau,\eta,\zeta) = 
\hat{\chi}^i\left(g_{i,\eta,\zeta}^{-1}(s),\eta,\zeta\right)$ and $\bar{\zeta}^i(\tau,\eta,\zeta) = 
\hat{\zeta}^i\left(g_{i,\eta,\zeta}^{-1}(s),\eta,\zeta\right)$,
\begin{align}
\int\limits_0^{s_F^i(\eta,\zeta)} \Phi(\hat{x}^i(s,\eta,\zeta), 
\hat{\xi}^i(s,\eta,\zeta),\eta,\zeta)^\ell 
ds & =
\nonumber 
\\
\int\limits_{\Phi(x,\xi,\eta,\zeta)}^{\Phi\left(x_F^i(\eta,\zeta),\xi_F^i(\eta,\zeta),\eta,\zeta\right)} 
\frac{\tau^\ell 
d\tau}{g(\bar{\chi}^i(\tau,\eta,\zeta), \bar{\zeta}^i(\tau,\eta,\zeta),\eta,\zeta)} & \leq \nonumber \\
\frac{1}{m_\mu e^\gamma - M_\mu 
e^{e^{-\gamma}-\gamma}}\frac{\Phi(x,\xi,\eta,\zeta)^{\ell+1}}{\ell+1},
\end{align}
and hence, \eqref{eq:MPhi} holds by the choice of $m_\Phi$ in \eqref{eq:mphi}.

The relations \eqref{eq:infsaest} now follow by using \eqref{eq:infsaestpw} in the successive 
approximations of $\Delta K_\ell^i$ and $\Delta L_\ell^i$ and taking the $L^2$ norm over 
$\eta,\zeta$. Due to linearity, the integral equations for $\Delta K_\ell^i$ and $\Delta L_\ell^i$, 
are 
of the same form as \eqref{eq:infkie}, and the choice of $M$ in \eqref{eq:seM} guarantees that  
\eqref{eq:infsaest} is satisfied for $\ell =0$. For any $\ell > 0$, we insert \eqref{eq:infsaestpw} 
to 
the integral equations for $\Delta K_\ell^i$ and $\Delta L_\ell^i$ and use \eqref{eq:MPhi} 
together 
with the choice of $M_{K,L}$ in \eqref{eq:seMKL} to show that  \eqref{eq:infsaestpw} holds for 
$\ell+1$. Finally, the relations \eqref{eq:infsaest} follow by taking the $L^2$ norm over 
$\eta,\zeta$ in \eqref{eq:infsaestpw}.
\end{proof}
\end{lemma}

\paragraph*{Proof of Theorem~\ref{thm:infkwp}} By Lemma~\ref{lem:sa}, the sequences of 
successive approximations for the kernels $K^i$ and $L^i$ converge in $\mathcal{H}_i$ for all 
$i \in \{a,b,c\}$ in the sense of \eqref{eq:infksac}, which shows the existence and well-posedness 
of the solutions $K^i, L^i$ to the kernel equations \eqref{eq:infkseg}--\eqref{eq:infkbcsegc}, 
which then uniquely determine the solution to the kernel equations \eqref{eq:infk}--\eqref{eq:infl} 
in the stated sense, i.e., $K, L \in L^\infty(\mathcal{T}; L^2([0,1]^2; \mathbb{R}))$.

\section{Numerical Examples and Simulations} \label{sec:numex}

\subsection{Illustration of Theorem~\ref{thm:stab} and Theorem~\ref{thm:infUappr1}} 
\label{sec:numex1}

\begin{figure}[!htb]
\begin{center}
\includegraphics[width=\columnwidth]{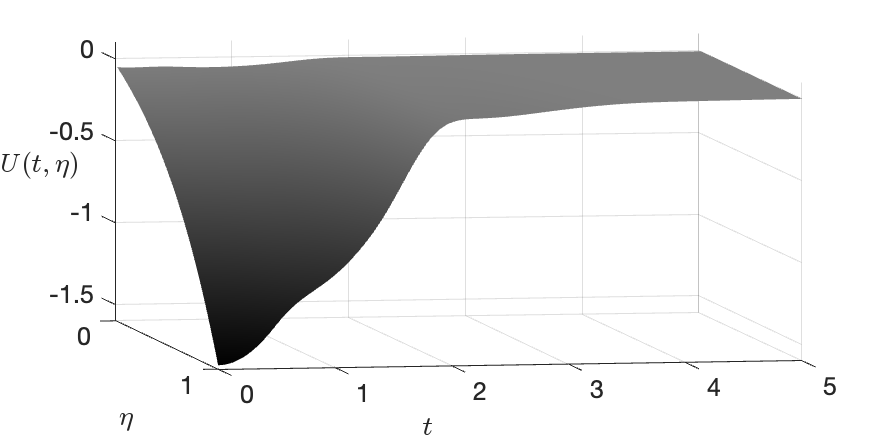}
\end{center}
\caption{The controls $U(t,\eta)$ from \eqref{eq:infU} for $t \in [0,5]$ and $\eta \in [0,1]$.}
\label{fig:exU}
\end{figure}

For a numerical example, consider the following parameters for $x,y,\eta,\zeta \in [0,1]$
\begin{subequations}
\label{eq:exparam}
\begin{align}
\lambda(x,y)& = 1, \qquad \mu(x,\eta) = 2-\eta, \\
\sigma(x,y,\zeta) & = W(x,y,\zeta) = (x+1)y\left(\zeta + \frac{1}{2}\right), \\
\theta(x,\eta,\zeta) & = \sigma(x,\eta,\zeta), \qquad 
\psi(x,\eta,\zeta) = \eta - \zeta, \\
Q(y,\zeta) & = \left(y + \frac{1}{2}\right)\zeta, \qquad R(\eta,\zeta) = 0.
\end{align}
\end{subequations}
For illustration of Theorem~\ref{thm:stab}, the continuum system \eqref{eq:inf}, \eqref{eq:infbc} 
with parameters 
\eqref{eq:exparam} is approximated by a grid of $50$ points in $y,\eta,\zeta$ and $128$ points 
in $x$, where we use finite differences to approximate the differential operators. The kernels $K, 
L$ for the control law \eqref{eq:infU} are approximated by 4-D power series of order ten by 
extending the power series approach from \cite{HumBek25, VazCheCDC23} to  4-D, and 
thereafter evaluating the obtained kernels at the employed grid points for computing the control 
law \eqref{eq:infU}. The initial conditions for the simulation are taken as $u_0(x,y) \equiv 
\int\limits_0^1 Q(y,\zeta)d\zeta$ and $v_0(x,\eta)\equiv 1$, and the closed-loop ODE resulting 
from the approximation is simulated using \texttt{ode45} in MATLAB. 
The control $U(t,\eta)$ based on \eqref{eq:infU} for $t \in [0,5]$ and $\eta \in [0,1]$ in the 
simulation is shown in Fig.~\ref{fig:exU}. One can see that the control input tends to zero 
exponentially and it is very close to zero by $t=5$ in the simulation. Since the control input 
contains a weighted average of the solution components, one can conclude that the closed-loop 
system is exponentially stable. We note that, based on numerical simulations, the open-loop 
system is unstable.

\begin{figure}[!htb]
\begin{center}
\includegraphics[width=\columnwidth]{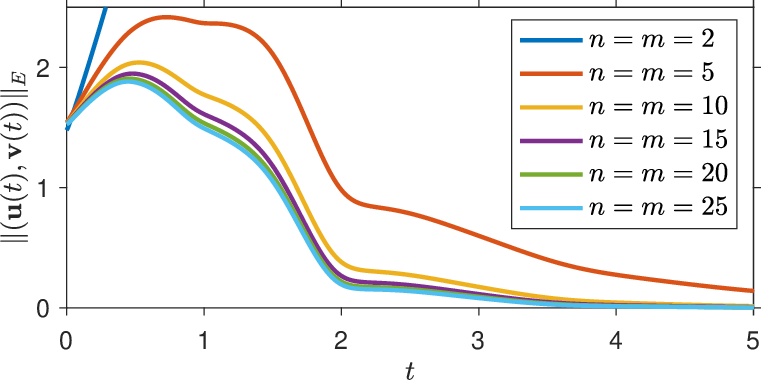}
\end{center}
\caption{Norm of the solution of the $n+m$ system for different $n=m$ under the 
continuum-kernels-based control law \eqref{eq:infUappr1} for $t \in [0, 5]$.}
\label{fig:nmcomp}
\end{figure}

In order to illustrate Theorem~\ref{thm:infUappr1}, we view the continuum parameters 
\eqref{eq:exparam} as continuum approximations of the respective $n+m$ parameters, defined 
for $i,\ell = 1,\ldots,n$ and $j,p = 1,\ldots,m$, as $r_{j,\ell} = 0$ and
\begin{subequations}
\label{eq:exnmparam}
\begin{align}
\lambda_i& = 1, \qquad \mu_j = 2-\frac{j}{m}, \\
\sigma_{i,\ell}(x) & = (x+1)\frac{i}{n}\left(\frac{\ell}{n} + \frac{1}{2}\right), \\
w_{i,p}(x) & = (x+1)\frac{i}{n}\left(\frac{p}{m} + \frac{1}{2}\right), \\
\theta_{j,\ell}(x) & =(x+1)\frac{j}{m}\left(\frac{\ell}{n} + \frac{1}{2}\right), \\
\psi_{j,p} & = \frac{j}{m} - \frac{p}{m}, \qquad q_{i,p} = \left(\frac{i}{n} + 
\frac{1}{2}\right)\frac{p}{m},
\end{align}
\end{subequations}
where we consider various $n,m$ to illustrate how they affect the closed-loop performance. We 
simulate the $n+m$ system with parameters 
\eqref{eq:exnmparam} for $n=m \in \{2,5,10,15,20,25\}$ under the continuum-kernels-based 
control 
law \eqref{eq:infUappr1}. The norm of the solution of the closed-loop system is displayed in 
Fig.~\ref{fig:nmcomp}, where one can see that the controller fails to stabilize the closed-loop 
system when $n=m=2$, and that when the closed-loop system is stable, the convergence rate is 
slower for smaller $n$ and 
$m$. This is expected, because the approximation accuracy of the continuum kernels is 
expected to deteriorate (when compared to the exact $n+m$ kernels) when $n$ and $m$ are 
small. We note here that as $n,m$ become larger, the $n+m$ kernels computation, 
based on the respective $n+m$ kernel equations from \cite{HuLDiM16}, may become intractable. 
This is because computing the exact $n+m$ kernels requires solving $m(n+m)$ (2-D) kernel 
equations, whereas computing the stabilizing, continuum-based kernels requires solving two 
(4-D) kernel equations, which is independent of $n$ and $m$.

As a heuristic approach to estimating 
the accuracy of the approximate control gains a priori, one can insert the sampled kernels 
\eqref{eq:infUappr1k} to the $n+m$ kernel 
equations \eqref{eq:nmk}, \eqref{eq:nmkbc} and test numerically, how accurately the equations 
are satisfied. We have
performed this for the case corresponding to the cases $n=m \in \{2,5,10\}$ in 
Fig.~\ref{fig:nmcomp}, which 
resulted in the residuals (maximum 
over $\mathcal{T}$ or a given boundary) shown in Table~\ref{tab:cc}. We note that the 
inconsistent residuals in \eqref{eq:nmk1}  (as the residual increases from the case $n=m=5$ to 
$n=m=10$) 
are caused by numerical inaccuracies related to the continuity condition 
\eqref{eq:infkbcsegc}. In more detail, this continuity condition is not satisfied exactly by the 
approximate (for the 
continuum) power 
series solution, which may (and in this case does) lead to non-smoothness in $K$ in the 
neighborhood of the characteristic hypersurface $\xi(\eta,\zeta) = \rho(x,\eta,\zeta)$. 
Consequently, such non-smoothness gets amplified by differentiation, which may (and in this 
case does) lead to large (local) residuals in \eqref{eq:nmk1} . Hence, the residuals shown in 
Table~\ref{tab:cc} may still be too conservative for accurately estimating approximation errors in 
the control gains, so that, e.g., the integral form of \eqref{eq:nmk}, \eqref{eq:nmkbc} should be 
considered instead for computing more reliable error estimates.

\begin{table}[!h] 
\begin{center}
\begin{tabular}{|c|cccccc|} \hline & & & & & &  \\ [-8pt]
 & \eqref{eq:nmk1} & \eqref{eq:nmk2} & \eqref{eq:nmkbc1} & \eqref{eq:nmkbc2}  & 
 \eqref{eq:nmkbc3} & \eqref{eq:nmkbc4} \\
\hline & & & & & &  \\ [-8pt]
$n=m=2$ & 1.01 & 0.99 & 0.37 & 0.30 & 0.37 & 0.14 \\
$n=m=5$ & 0.97 & 0.50 & 0.19 & 0.31 & 0.10 & 0.11 \\
$n=m=10$ & 1.11 & 0.40 & 0.14 & 0.31 & 0.05 & 0.09 \\
\hline 
\end{tabular}
\end{center}
\caption{Residuals of the $n+m$ kernel equations \eqref{eq:nmk}, \eqref{eq:nmkbc} for the 
kernels employed in 
Fig.~\ref{fig:nmcomp}.}
\label{tab:cc}
\end{table}

\subsection{Illustration of Proposition~\ref{thm:ave}}

For illustrating  Proposition~\ref{thm:ave}, we consider an 
$n+m$ system with parameters, 
defined for $i,\ell = 1,\ldots,n$ and $j,p = 1,\ldots,m$, as $r_{j,\ell} = 0,  q_{i,p} = 1,$ and
\begin{subequations}
\label{eq:exnmparam2}
\begin{align}
\lambda_i& = 1, \qquad \mu_j = 1-\frac{j}{2m}, \\
\sigma_{i,\ell}(x) & = x\frac{i}{n}\left(\frac{\ell}{n} + \frac{1}{2}\right), \quad
w_{i,p}(x)  = x\frac{i}{n}\left(\frac{p}{m} + \frac{1}{2}\right), \\
\theta_{j,\ell}(x) & =x\frac{j}{m}\left(\frac{\ell}{n} + \frac{1}{2}\right), \quad
\psi_{j,p} = \frac{j}{2m} - \frac{p}{2m}
\end{align}
\end{subequations}
where we take $n=m=10$. We construct two different approximations for \eqref{eq:exnmparam2} 
and the respective controllers \eqref{eq:infUappr2} to illustrate 
Proposition~\ref{thm:ave}. Firstly, 
we construct a continuum approximation of the $n+m$ system with parameters 
\eqref{eq:exnmparam2}, using Remark~\ref{rem:nmcont} and choosing, for $x,y,\eta,\zeta \in 
[0,1]$, the continuum parameters as
\begin{subequations}
\label{eq:exnmparam2c}
\begin{align}
\lambda(x,y) & = 1, \qquad \mu(x,\eta) = 1-\frac{1}{2}\eta, \\
\sigma(x,y,\zeta) & = W(x,y,\zeta) = xy\left(\eta + \frac{1}{2}\right), \\
\theta(x,\eta,\zeta) & = \sigma(x,\eta,\zeta), \qquad \psi(x,\eta,\zeta) = \frac{1}{2}(\eta- \zeta), \\
Q(y,\zeta) & = 1, \qquad R(\eta,\zeta) = 0.
\end{align}
\end{subequations}
Respectively, the continuum kernel equations are solved similarly to Section~\ref{sec:numex1}. 
Secondly, we construct an average approximation with states $\bar{u}, \bar{v}$, and parameters 
$\bar{r} = 0, \bar{q} = 1, \bar{\psi} = 0$, and
\begin{subequations}
\label{eq:exnmparam2ave}
\begin{align}
\bar{\lambda} & = 1, \qquad \bar{\mu}= \frac{3}{4}, \\
\bar{\sigma}(x) & = \bar{W}(x) = \bar{\theta}(x) = \frac{1}{2}x,
\end{align}
\end{subequations}
which are obtained by taking the averages of the respective $n+m$ parameters over 
$i,\ell,j,p$. The respective $1+1$ (continuum) kernels \eqref{eq:2x2k}, \eqref{eq:2x2kbc} are 
solved using finite differences and successive approximations.

\begin{figure}[!htb]
\begin{center}
\includegraphics[width=\columnwidth]{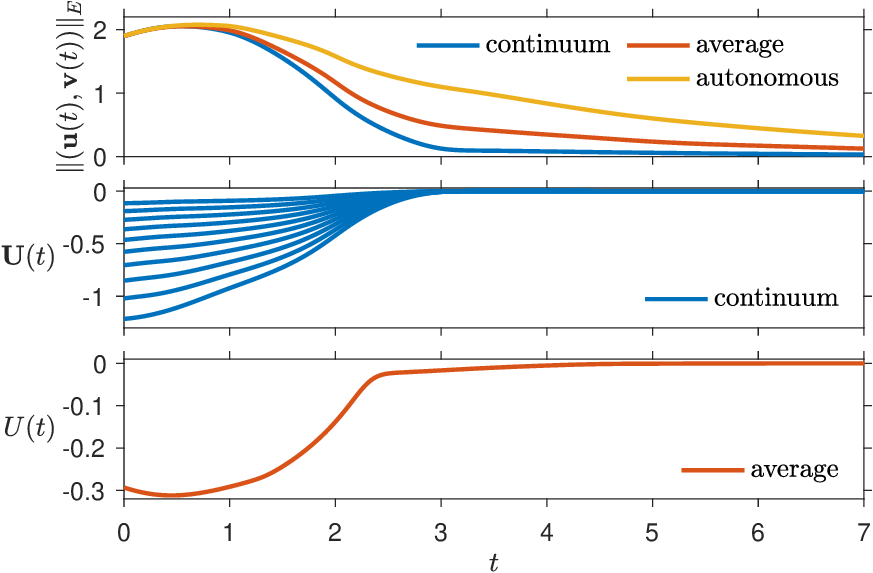}
\end{center}
\caption{Uppermost plot: norm of the solution of the $n+m$ system under the
 control law \eqref{eq:infUappr2}, employing measurements/kernels, according to 
 \eqref{eq:infUappr2k} with \eqref{eq:uvave}, and 
\eqref{eq:2x2k}, \eqref{eq:2x2kbc} with $\tilde{u} \equiv \bar{u}, \tilde{v} \equiv \bar{v}$, from the 
continuum and 
 average approximations, respectively, and a  comparison to the autonomous system's solution. 
 Lower  plots: the respective control inputs \eqref{eq:infUappr2} (where each component of 
 $\mathbf{U}$ coincides with $U$ in the bottom plot).}
\label{fig:ave43}
\end{figure}

For the simulations, the initial conditions are taken as $u_0^i = 0.9$ for $i = 1,\ldots,n$ and 
$v_0^j  =1$ for $j = 1,\ldots,m$ for the $n+m$ system and as $u_0 = v_0 = 1 = \bar{u}_0 = 
\bar{v}_0$ for the continuum approximations. The simulation results are shown in the uppermost 
plot of Fig.~\ref{fig:ave43}, where the norm of the solution of the $n+m$ system under the control 
law \eqref{eq:infUappr2}, which employs kernels/measurements from each of the two constructed 
continuum approximations (according to \eqref{eq:infUappr2k} with \eqref{eq:uvave}, and 
\eqref{eq:2x2k}, \eqref{eq:2x2kbc} with $\tilde{u} \equiv \bar{u}, \tilde{v} \equiv \bar{v}$, 
respectively), is compared 
with the norm of the solution of the autonomous $n+m$ system. One can see that the controls, 
which are shown in the lower plots of Fig.~\ref{fig:ave43}, improve the transient response of 
the $n+m$ system by improving the convergence rate of the solution as compared to the solution 
of the autonomous system. 
Moreover, this improvement is more evident under the continuum approximation-based controls, 
which is expected, because the continuum \eqref{eq:exnmparam2c} provides a better 
approximation of the $n+m$ parameters \eqref{eq:exnmparam2} than the (very simple) average 
system with parameters \eqref{eq:exnmparam2ave}. This improvement would be even more 
pronounced under the continuum approximation-based controls as $n$ and $m$ increase, since 
the respective solutions' approximation accuracy improves. We note that the autonomous $n+m$ 
system with parameters \eqref{eq:exnmparam2} is exponentially stable in the simulation, so that 
the solution tends eventually to zero even in the absence of controls, albeit the decay rate is quite 
small. 

\section{Conclusions and Discussion} \label{sec:conc}

The paper considered different micro-macro control scenarios for large-scale $n+m$ and 
continuum $\infty+\infty$ hyperbolic systems. Firstly, we derived in a constructive manner the 
class of $\infty+\infty$ 
hyperbolic PDEs as a continuum approximation of large-scale $n+m$ hyperbolic PDEs, and then 
solved the backstepping state-feedback stabilization problem for the
$\infty+\infty$ PDEs. In particular, we established well-posedness of the resulting 4-D continuum 
kernel equations and closed-loop stability constructing a Lyapunov functional. Secondly, we 
developed micro-macro controllers for large-scale $n+m$ 
systems based on control kernels and/or measurements  obtained on the basis of 
the $\infty+\infty$ continuum system. In particular, we established that the macro 
measurements/kernels can approximate the micro measurements/kernels in certain sense, which 
then allowed us to derive specific stability properties for the respective closed-loop systems 
utilizing infinite-dimensional ISS arguments. The effectiveness of the proposed controllers was 
illustrated in numerical simulations.

Among the different research problems one can study capitalizing on the results of the present 
paper, we discuss the following two. The first is development of systematic computational tools 
for solving the continuum, 4-D kernel equations, in order to maximize the potential benefits in 
computational complexity of computation of stabilizing kernels for large-scale and continua-of 
hyperbolic systems. For example, in our numerical example in Section~\ref{sec:numex1}, we 
solved the 
continuum kernel equations using 4-D power series (inspired from \cite{VazCheCDC23, 
HumBek25}). However, such a 
practical approach may not be the optimal choice, as, even though the respective computational 
complexity does not scale with $n$ and $m$, it still grows with the order $N$ of the power series 
needed to guarantee stabilization, as $\mathcal{O}(N^4)$. (This problem is also related to the 
study of optimal construction of the continuum approximation.) The second is the application of 
the 
control design methods 
developed here to specific engineering applications, in particular, to lane-free (or multi-lane) and 
continuum/multi-class traffic (see \cite{MalPap21, YuHKrs21, BurYuH21, DeoSch09}). In particular, 
it is anticipated that lane-free and continuum-class traffic flow models to be possible, in principle, 
to be recast in the form of the continuum systems considered in the present paper.

\appendix

\renewcommand{\theequation}{\thesubsection.\arabic{equation}}
\renewcommand{\thetheorem}{\thesubsection.\arabic{theorem}}

\subsection{Derivation of Kernel Equations} \label{app:ker}

\setcounter{equation}{0}
\setcounter{theorem}{0}

Let us first differentiate \eqref{eq:infV2} with respect to $x$ and use the Leibniz
rule to get
\begin{align}
	\beta_x(t,x,\eta) & = v_x(t,x,\eta) -
	\int\limits_0^1L(x,x,\eta,\zeta)v(t,x,\zeta)d\zeta
	\nonumber \\
	& \qquad - \int\limits_0^1K(x,x,\eta,\zeta)u(t,x,\zeta)d\zeta \nonumber \\
	& \qquad  -
	\int\limits_0^x\int\limits_0^1L_x(x,\xi,\eta,\zeta)v(t,\xi,\zeta)d\zeta d\xi
	\nonumber \\
	& \qquad  - \int\limits_0^x \int\limits_0^1K_x(x,\xi,\eta,\zeta)u(t,\xi,\zeta)d\zeta d\xi.
\end{align}
Moreover, differentiating \eqref{eq:infV2} with respect to $t$
and using \eqref{eq:inf2} gives
\begin{align}
	\beta_t(t,x,\eta)
	& = \mu(x,\eta)v_x(t,x,\eta) +
	\int\limits_0^1\theta(x,\eta,\zeta)u(t,x,\zeta)d\zeta 
	\nonumber \\
	& \quad+ \int\limits_0^1\psi(x,\eta,\zeta)v(t,x,\zeta)d\zeta  \nonumber \\
	& \quad - \int\limits_0^x\int\limits_0^1L(x,\xi,\eta,\zeta)\mu(\xi,\zeta)v_\xi(t,\xi,\zeta)d\zeta 
	d\xi
	\nonumber \\
	& \quad - \resizebox{.72\columnwidth}{!}{$\displaystyle \int\limits_0^x\int\limits_0^1
	L(x,\xi,\eta,\zeta)\int\limits_0^1\theta(\xi,\zeta,\chi)u(t,\xi,\chi)d\chi d\zeta d\xi$}
	\nonumber \\
	& \quad
	-\resizebox{.72\columnwidth}{!}{$\displaystyle 
	\int\limits_0^x\int\limits_0^1L(x,\xi,\eta,\zeta)\int\limits_0^1\psi(\xi,\zeta,\chi)v(t,\xi,\chi)d\chi
	d\zeta d\xi$}
	\nonumber \\
	& \quad + \int\limits_0^x
	\int\limits_0^1K(x,\xi,\eta,\zeta)\lambda(\xi,\zeta)u_{\xi}(t,\xi,\zeta)d\zeta
	d\xi \nonumber \\
	& \quad -
	\resizebox{.72\columnwidth}{!}{$\displaystyle \int\limits_0^x\int\limits_0^1K(x,\xi,\eta,\zeta)
		\int\limits_0^1\sigma(\xi,\zeta,\chi)u(t,\xi,\chi)d\chi d\zeta d\xi$}
	\nonumber \\
	& \quad -
	\resizebox{.72\columnwidth}{!}{$\displaystyle\int\limits_0^x\int\limits_0^1K(x,\xi,\eta,\zeta)\int\limits_0^1W(\xi,\zeta,\chi)v(t,\xi,\chi)d\chi
	 d\zeta	d\xi$},
\end{align}
where integration by parts further gives 
\begin{align}
	\int\limits_0^xL(x,\xi,\eta,\zeta)\mu(\xi,\zeta)v_\xi(t,\xi,\zeta)d\xi
	& = \nonumber \\
	\resizebox{.92\columnwidth}{!}{$\displaystyle L(x,x,\eta,\zeta)\mu(x,\zeta)v(t,x,\zeta)
	- L(x,0,\eta,\zeta)\mu(0,\zeta)v(t,0,\zeta)$} & \nonumber \\
\resizebox{.92\columnwidth}{!}{$\displaystyle	- \int\limits_0^x 
\left(L_{\xi}(x,\xi,\eta,\zeta)\mu(\xi,\zeta) +
	L(x,\xi,\eta,\zeta)\mu_{\xi}(\xi,\zeta)\right)v(t,\xi,\zeta)d\xi$},
\end{align}
and
\begin{align}
	\int\limits_0^x K(x,\xi,\eta,\zeta)\lambda(\xi,\zeta)u_{\xi}(t,\xi,\eta,\zeta)
	d\xi
	& = \nonumber \\
	\resizebox{.92\columnwidth}{!}{$\displaystyle K(x,x,\eta,\zeta)\lambda(x,\zeta)u(t,x,\zeta) -
	K(x,0,\eta,\zeta)\lambda(0,\zeta)u(t,0,\zeta)$}
	& \nonumber \\
	\resizebox{.92\columnwidth}{!}{$\displaystyle- 
	\int\limits_0^x(K_{\xi}(x,\xi,\eta,\zeta)\lambda(\xi,\zeta) +
	K(x,\xi,\eta,\zeta)\lambda_{\xi}(\xi,\zeta))u(t,\xi,\zeta)d\xi$}.
\end{align}
Hence, we  get kernel equations \eqref{eq:infk} with boundary conditions \eqref{eq:infkbc1}, 
\eqref{eq:infkbc2}, and
\begin{subequations}
\begin{align}
\psi(x,\eta,\zeta) - L(x,x,\eta,\zeta)\mu(x,\zeta) + \mu(x,\eta)L(x,x,\eta,\zeta) & =  0, \\
\theta(x,\eta,\zeta) + K(x,x,\eta,\zeta)\lambda(x,\zeta) +  \mu(x,\eta)K(x,x,\eta,\zeta)  & = 0, \\
\int\limits_0^1 K(x,0,\eta,\zeta)\lambda(0,\zeta)\int\limits_0^1Q(\zeta,\chi)h(\chi)d\chi d\zeta 
& = \nonumber \\ \int\limits_0^1 L(x,0,\eta,\zeta)\mu(0,\zeta)h(\zeta)d\zeta - \int\limits_0^\eta 
G(x,\eta,\zeta)h(\zeta)d\zeta, 
\end{align}	
\end{subequations}
for all $h \in L^2([0,1]; \mathbb{R})$, where changing the order of integration and splitting the 
integrals over $\zeta \in [0,1]$ into $\zeta \in [0, \eta]$ and $\zeta \in 
(\eta,1]$ gives \eqref{eq:infkbc3}, and that $G$ is given, for $\zeta < \eta$, by
\begin{align}
\label{eq:infGdef}	
G(x,\eta,\zeta) & = L(x,0,\eta,\zeta)\mu(0,\zeta) \nonumber \\
& \qquad - \int\limits_0^1K(x,0,\eta,\chi)\lambda(0,\chi)Q(\chi,\zeta)d\chi.
\end{align}
Finally, inserting \eqref{eq:infV} to \eqref{eq:infts1} gives that $C^-$ and $C^+$ need to satisfy
\begin{subequations}
\label{eq:infC}	
\begin{align}
C^-(x,\xi,y,\zeta) & = \int\limits_0^1W(x,y,s)L(x,\xi,s,\zeta)ds \nonumber \\
& \qquad  + \int\limits_\xi^x \int\limits_0^1 
C^-(x,\chi,y,s)L(\chi,\xi,s,\zeta)dsd\chi, \label{eq:infC-} \\
C^+(x,\xi,y,\zeta) & = \int\limits_0^1W(x,y,s)K(x,\xi,s,\zeta)ds \nonumber \\
& \qquad + \int\limits_\xi^x\int\limits_0^1 C^-(x,\chi,y,s)K(\chi,\xi,s,\zeta)dsd\chi,
\label{eq:infC+}
\end{align}
\end{subequations}
where $C^+,C^- \in L^\infty(\mathcal{T}; L^2([0,1]^2; \mathbb{R}))$.

\begin{remark}
\label{rem:infk}
Analogously to the case of finite $m$ (see, e.g., \cite[Thm A.1]{HuLVaz19}, \cite[Lem. 
1]{HumBek26}), the boundary conditions on $(x,\xi) = (0,0)$ are (generally) 
over-determined (for $L$ on $\eta \leq \zeta$) because of \eqref{eq:infkbc1} and 
\eqref{eq:infkbc3}, 
\eqref{eq:infkbc2}, which stems a potential discontinuity in the $L$ kernels. Hence, the kernel 
equations \eqref{eq:infk}--\eqref{eq:infl} are given for almost all $(x,\xi) \in \mathcal{T}$ and 
$\eta,\zeta \in [0,1]$, 
such that $K,L \in L^\infty(\mathcal{T}; L^2([0,1]^2; \mathbb{R}))$\footnote{Assuming 
\eqref{eq:mupsiass} guarantees that \eqref{eq:infkbc1} is well-posed in this sense.}. In order to 
gain more regularity, the kernels can be segmented into subdomains \eqref{eq:Hseg}, according 
to the characteristic hypersurface of \eqref{eq:infk2}. The resulting segmented kernels are then 
continuous in $(x,\xi)$ and they satisfy the respective segmented kernel equations 
\eqref{eq:infkseg}--\eqref{eq:infkbcsegc}, where we have an additional continuity condition 
\eqref{eq:infkbcsegc} for the $K$ kernel, whereas the $L$ kernel (generally) has a discontinuity 
along its characteristic hypersurface.
\end{remark}

\subsection{Invertibility of the Backstepping Transformation \eqref{eq:infV}}

\setcounter{equation}{0}
\setcounter{theorem}{0}

\begin{lemma}
\label{lem:infkinv}		
Under Assumption~\ref{ass:inf}, the transformation \eqref{eq:infV} is boundedly invertible on 
$E_c$.

\begin{proof}
Consider an arbitrary, fixed $t \geq 0$, so that $\alpha(t), \beta(t), u(t), v(t) \in L^2([0,1];L^2([0,1]; 
\mathbb{R}))$ and $K, L \in L^\infty(\mathcal{T}; L^2([0,1]^2; \mathbb{R}))$. Inserting $u(t) = 
\alpha(t)$ from \eqref{eq:infV1} to \eqref{eq:infV2}, it remains to solve $v(t)$ from
\begin{align}
	\label{eq:Ttmp}
v(t,x,\eta) & =  \int\limits_0^x\int\limits_0^1 L(x,\xi,\eta,\zeta) v(t,\xi,\zeta)	d\zeta d\xi \nonumber \\
& \quad + \beta(t,x,\eta) + \int\limits_0^x\int\limits_0^1 K(x,\xi,\eta,\zeta) \alpha(t,\xi,\zeta)d\zeta 
d\xi 
\nonumber \\
& =: \mathcal{V}v(t,x,\eta).
\end{align}
Using similar arguments to the proof of \cite[Thm 2.3.5]{HocBook}, we show that there exists 
some $\ell > 0$ such that operator $\mathcal{V}^\ell$, where $\mathcal{V}$ is defined in 
\eqref{eq:Ttmp}, is a 
contraction on $L^2([0,1]; L^2([0,1]; \mathbb{R}))$. Let us
denote, for almost all $(x,\xi) \in \mathcal{T}$,
\begin{align}
L_1(x,\xi,\cdot) & = \int\limits_0^1 L(x,\xi,\cdot,\zeta)d\zeta, \\
L_\ell(x,\xi,\cdot) & = \int\limits_\xi^x \int\limits_0^1 L(x,s,\cdot,\zeta)L_{\ell-1}(s,\xi,\zeta)d\zeta ds,
\end{align}
where $\ell \geq 2$, so that
\begin{align}
\mathcal{V}^\ell v_1(t,x,\eta) - \mathcal{V}^\ell v_2(t,x,\eta) & =  \nonumber \\
\resizebox{!}{.038\textwidth}{$\displaystyle \int\limits_0^x \int\limits_0^1
L(x,\xi,\eta,\zeta)L_{\ell-1}(x,\xi,\zeta)\left(v_1(t,\xi,\zeta) - v_2(t,\xi,\zeta)\right) 
d\zeta d\xi$},
\end{align}
which holds in the $L^2([0,1]; L^2([0,1]; \mathbb{R}))$ sense  in terms of $(x,\eta)$. 
Now, let
\begin{align}
	M_{L_1} = \esssup_{(x,\xi) \in \mathcal{T}} \left\|\int\limits_0^1 L(x,\xi,\cdot,\zeta)d\zeta 
	\right\|_{L^2},
\end{align}
so that $\|L_1(x,\xi,\cdot)\|_{L^2} \leq M_{L_1}$ holds by construction, and let us make the 
induction assumption that 
\begin{align}
\|L_\ell(x,\xi,\cdot)\|_{L^2} \leq \frac{M_{L_1}^\ell(x-\xi)^{\ell-1}}{(\ell-1)!},
\end{align}
holds for some $\ell \in \mathbb{N}$, for almost all $(x,\xi) \in \mathcal{T}$. Now, by 
Cauchy-Schwarz inequality,
\begin{align}
\|L_{\ell+1}(x,\xi,\cdot)\|_{L^2} & \leq \int_\xi^x M_{L_1}\frac{M_{L_1}^\ell(s-\xi)^{\ell-1}}{(\ell-1)!}ds 
\nonumber \\
& = \frac{M_{L_1}^{\ell+1}(x-\xi)^{\ell}}{\ell!},
\end{align}
and hence,
\begin{align}
\label{eq:Tpf}
\|\mathcal{V}^\ell v_1(t) - \mathcal{V}^\ell v_2(t)\|_{E_c} & \leq \nonumber \\
\frac{M_{L_1}^\ell}{(\ell-1)!} \left\| \int\limits_0^x\int\limits_0^1 \left(
v_1(t,\xi,\zeta) - v_2(t,\xi,\zeta)\right)d\zeta d\xi \right\|_{L^2} & \leq \nonumber \\
\frac{M_{L_1}^\ell}{(\ell-1)!} \|v_1(t) - v_2(t)\|_{E_c}, &
\end{align}
where $\displaystyle \frac{M_{L_1}^\ell}{(\ell-1)!} < 1$ for sufficiently large $\ell$, so that 
$\mathcal{V}^\ell$ is a contraction for any such $\ell$. Hence, \eqref{eq:Ttmp} has a unique 
solution in $L^2([0,1]; L^2([0,1]; \mathbb{R}))$ by \cite[Thm 2.1.2]{HocBook}.
\end{proof}
\end{lemma}

\end{document}